\documentclass[12pt,reqno]{amsart}
\usepackage{amsmath,amsfonts,amsthm,amssymb,color}
\usepackage[T1]{fontenc}
\usepackage{enumitem}
\usepackage{hyperref}
\hypersetup{
     colorlinks   = true,
     citecolor    = blue,
     linkcolor    = blue
}
\usepackage[left=1.1in, right=1.1in, top=1.1in, bottom=1.1in]{geometry}
\usepackage[font=small]{caption}

\usepackage{pdfsync}

\usepackage{amsmath}
\usepackage{amssymb}
\usepackage{fancyhdr}
\usepackage{graphicx}
\usepackage{mathrsfs}
\usepackage{amsthm}
\usepackage{cleveref}

\usepackage{xcolor}
\usepackage{lastpage}
\usepackage{fancyvrb}
\usepackage{listings}
\usepackage{bbm}

\usepackage{soul}
\usepackage{mathtools}

\usepackage{tikz}

\newcommand{\dd}{\mathrm{d}}

\newcommand{\iot}{\int_{0}^{t}}
\newcommand{\ios}{\int_{0}^{s}}

\newcommand{\norm}[1]{\left\lVert#1\right\rVert}

\newcommand{\cn}{\mathcal N}

\newcommand{\E}{\mathbb{E}}
\newcommand{\PP}{\mathbb{P}}
\newcommand{\R}{\mathbb{R}}

\newcommand{\Pbf}{\mathbf{P}}
\newcommand{\Pcal}{\mathcal{P}}
\newcommand{\Ebf}{\mathbf{E}}

\newcommand{\al}{\alpha}
\newcommand{\ep}{\varepsilon}

\newcommand{\ga}{\gamma}
\newcommand{\ka}{\kappa}
\newcommand{\la}{\lambda}

\newcommand{\lc}{\left[}
\newcommand{\rc}{\right]}

\newcommand{\cox}{\text{Cox}}

\newtheorem{theorem}{Theorem}[section]

\newtheorem{definition}[theorem]{Definition}

\newtheorem{hypothesis}[theorem]{Hypothesis}
\newtheorem{lemma}[theorem]{Lemma}
\newtheorem{notation}[theorem]{Notation}

\newtheorem{proposition}[theorem]{Proposition}

\theoremstyle{remark}
\newtheorem{remark}[theorem]{Remark}

\crefname{appendix}{Appendix}{Appendices}
\crefname{corollary}{Corollary}{Corollaries}
\crefname{definition}{Definition}{Definitions}
\crefname{figure}{Figure}{Figures}
\crefname{hypothesis}{Hypothesis}{Hypotheses}
\crefname{lemma}{Lemma}{Lemmas}
\crefname{notation}{Notation}{Notations}
\crefname{proposition}{Proposition}{Propositions}
\crefname{remark}{Remark}{Remarks}
\crefname{section}{Section}{Sections}
\crefname{theorem}{Theorem}{Theorems}

\title[Diffusion Approximations for $\cox/G_t/\infty$ Queues]{Diffusion Approximations for $\cox/G_t/\infty$ Queues In A Fast Oscillatory Random Environment}

\author[Y. Liu \and H. Honnappa \and S. Tindel \and N. K. Yip]
{Yiran Liu \and Harsha Honnappa \and Samy Tindel \and  Nung Kwan Yip}

\address{Harsha Honnappa: School of Industrial Engineering, Purdue
University, 315 N. Grant Street, W. Lafayette, IN 47907, USA.}

 \address{Yiran Liu, Samy Tindel, Nung Kwan Yip: Department of Mathematics,
Purdue University,
150 N. University Street,
W. Lafayette, IN 47907,
USA.}
\keywords{Infinite server queue, random environment, homogenization.}

\email{[honnappa, liu387, stindel, yipn]@purdue.edu}

\thanks{S. Tindel is supported by the NSF grant
  DMS-1952966. H. Honnappa is supported by the NSF grant CMMI-1636069.}
  
 \subjclass[2010]{60K25, 60G55, 60F15, 90B22}

\begin{document}

\begin{abstract}
We study infinite server queues driven by Cox processes in a fast oscillatory random environment. While exact performance analysis is difficult, we establish diffusion approximations to the (re-scaled) number-in-system process by proving functional central limit theorems (FCLTs) using a stochastic homogenization framework. This framework permits the establishment of quenched and annealed limits in a unified manner. 
At the quantitative level, we identity two parameter regimes, 
termed subcritical and supercritical indicating the relative dominance between the
two underlying stochasticities driving our system:
the randomness in the arrival intensity and that in the serivce times.
We show that while quenched FCLTs can only be established in the subcritical 
regime, annealed FCLTs can be proved in both cases.
Furthermore, the limiting diffusions in the annealed FCLTs display qualitatively different diffusivity properties in the two regimes, even though the stochastic primitives are identical. In particular, when the service time distribution
is heavy-tailed, the diffusion is sub- and super-diffusive in the 
sub- and super-critical cases. The results illustrate intricate interactions 
between the underlying driving forces of our system.
\end{abstract}

\maketitle

\section{Introduction}

Nonstationary queueing models have been extensively studied, with particular focus on the Markovian setting. A typical assumption in this setting is that the arrival and service intensities are smooth, deterministic time-varying functions~\cite{whitt1}. However, in practice, queueing systems are often subject to `environmental' noise: for instance, while arrival intensities to call centers and hospitals display time-of-day (or `diurnal') effects, the intensity functions also vary based on the day-of-week and seasonal effects. In other queueing systems, particularly those with high intensity arrivals such as computer networks or cloud service systems, there is also intra-day and intra-hour stochastic variation in the intensity process. The performance of these queueing systems is therefore affected by both the smaller time-scale stochastic variations, as well as (relatively) longer time-scale time-of-day effects. In these settings, it is natural to model the queue as operating in a random environment.

 In this paper, we focus on the case where the intensity of the traffic process is a stochastic process itself. Specifically, we assume that the traffic process is a doubly stochastic Poisson process or Cox process. While much is known about nonstationary infinite server queues, infinite server queues driven by Cox processes are harder to analyze. 
In our previous work~\cite{HLTY}, a stochastic fluid approximation to the number-in-system process of a $\text{Cox}/G_t/\infty$ queue was established in the setting where the random environment is `fast oscillatory' or `rapidly averaging'.
More precisely, the number-in-system process was shown to converge to 
that of a $M_t/G_t/\infty$ queue. 
The results were obtained as a consequence of a stochastic averaging principle satisfied by the random environment. 
In the current paper, we establish stochastic diffusion 
approximations (in the form of functional central limit theorems) for a 
{\em rescaled} number-in-system process. We incorporate various
ingredients not considered in \cite{HLTY}. Most importantly, we consider
the speeding-up of the arrival intensity. 
Because of this, the system operates
in vastly different regimes depending on the intricate interaction between the
degree of speed-up and the tail behavior of the service time distribution. We still work under the assumption of separation of time scales
so that stochastic averaging principle for the random environment is 
applicable.

 Diffusion approximations of $\cox/G_t/\infty$ queues can be established using the two-parameter heavy-traffic limits developed in~\cite{PW} which 
considers a sequence of infinite server queues with increasing mean arrival rate. In other words, the setting assumes a sequence of stochastic intensity processes with increasing mean `level', engendering a notion of `heavy-traffic' in this infinite server setting. That work shows that the diffusion approximation to the re-scaled number-in-system process is a non-Markov Gaussian process. 
In order to prove this result, they augment the state space by tracking at each time $t$, 
the number of active jobs in the system with residual service time bigger
than certain level $y$. When applied to $\cox/G_t/\infty$ queues, this yields a heavy-traffic diffusion approximation under the so-called {\it annealed measure} which averages over the measure corresponding to the stochastic intensity process). See also~\cite{koops2017networks} for an application to shot noise driven $\cox/G/\infty$ queues. 
 
Indeed, the annealed measure is the default analytical framework for proving fluid and diffusive approximations of  $\cox/G/\infty$ queueing systems.
From a performance analysis perspective, this is an apparently sensible choice owing to the fact that the stochastic intensity process is latent (or unobserved) and the {specific} sample path of the intensity process to be experienced is unknown {\it a priori}. On the other hand, this also implies that the asymptotic approximations truly hold under the expectation that the queueing system will encounter {`typical'} environmental conditions. In practice, this is rarely the case. Alternatively, {conditioning} on the specific sample path of the random environment would yield asymptotic approximations under the so-called {\em quenched measure}, that hold {\it for (almost surely) any realized sample path of the random environment}. This would mean that
quenched analysis is impossible under the two-parameter heavy-traffic setting in~\cite{PW} 
wherein a sequence of doubly stochastic arrival processes is considered. Instead, the approximations must be demonstrated for a fixed sample path of the random environment. This raises the question of how to rescale the $\cox/G_t/\infty$ queue and establish asymptotic approximations under both the quenched and annealed measures in a unified manner. 

This paper addresses these desiderata by using a stochastic homogenization or averaging approach when the random environment is modeled by a general ergodic stochastic process. Herein the time-scale of the random environment model is `accelerated'  in relation to the natural time-scale of the queue so that the fluctuations of the random environment are averaged out. 
The specific contributions of the current paper include the following:

\begin{enumerate}[leftmargin=*]
 	\item 
We identify two parameter regimes, termed
{\it subcritical \emph{and} supercritical}, which 
relate the spatial scaling exponent $\beta$ for the (rescaled) stochastic intensity of the $\cox$ process and the (strong) ergodicity rate of the random 
environment process driving the intensity. 
Specifically, the subcritical and supercritical regimes correspond 
to $\beta < 1$ and $> 1$, respectively. 
(The boundary case $\beta=1$ is called critical.)
For both subcritical and supercritical cases, we can establish functional 
law of large numbers (FLLN) and functional central limit theorem (FCLT).

 	\item For the subcritical regime, 
we prove an FCLT in the quenched setting establishing the diffusion approximation for the number-in-system (or `state') process of the $\cox/G_t/\infty$ queue.
The result shows that, for almost every sample path of the random environment, the rescaled state process converges to a {\em subdiffusive} 
Gaussian process 
limit with fractional Brownian motion (FBM)-like structure when the service time distribution is `heavy-tailed'. On the other hand, it converges to a stationary ergodic diffusion type limit when the service time distribution is `light-tailed'. 
We also show that analysis in the quenched setting is not possible in the 
supercritical case. 

 	\item For both subcritical and supercritical regimes,
we also prove FCLTs in the annealed setting. Interestingly, the limit process in the subcritical case coincides with that in the quenched setting. This follows from the fact that the covariance structure of the Gaussian limit in the quenched setting does {\it not} depend on the specific sample path of the random environment. In other words, the approximation is robust to the specific sampled random environment. 
On the other hand, in the supercritical case,
the limiting stochastic process is {\em superdiffusive} 
when the service times are `heavy-tailed'. 
This suggests that the random environment can induce larger fluctuations.
This phenomenon has not been noted before in the literature on queues in random environments and is a clear manifestation of the intricate 
interaction between arrival intensity and service time distributions.
 \end{enumerate}
 
The results in this paper, therefore, align along three axes: 
(i) scaling regimes (sub- vs supercritical), 
(ii) random environment conditions (quenched vs annealed), and 
(iii) service time distributions (heavy- vs light-tailed). 
Table~\ref{tab:summary} below consolidates the FCLT results in this paper. As noted above, in the supercritical regime a quenched FCLT does not exist. The table exposes the discovery that with heavy-tailed (HT) service times (in a sense made clear following Hypothesis~\ref{hyp:F-bar-increments}), 
the FCLT yields diffusion approximations that are qualitatively different in the sub- and supercritical regimes. This can potentially have important 
implications for performance analysis and control of queues. However, 
we do not expand on these issues in this paper. Also note that we have not treated the critical (boundary) case ($\beta=1$) in this article for the sake of conciseness. Our conjecture is that it would yield both quenched and annealed CLTs, but with different limits. 
We note that the situation in Table~\ref{tab:summary} is reminiscent of 
\cite{Peterson08} which investigates random walks in random environments 
(RWREs). In Sections 2.2 and 2.3 of \cite{Peterson08}, different regimes of 
CLTs 
are evidenced according to a parameter $s$ of the environment distribution.
In that work, the critical case (corresponding to $s=2$) is also not
handled.
See however \cite[Theorem 3]{DGQuenched} for a description
of a quenched statement in the critical case.

\begin{table}[h]
\centering
\begin{tabular}{| c |c | c |}
\hline
Regimes  & Quenched & Annealed \\
  \hline
  \begin{tabular}{c}
Subcritical \\
$(\beta < 1)$
\end{tabular}
 & \begin{tabular} {c}HT : Subdiffusive\\ LT : Normal \end{tabular}  & \begin{tabular}{c} HT : Subdiffusive\\ LT : Normal \end{tabular}  \\
\hline
  \begin{tabular}{c}
Supercritical \\
$(\beta > 1)$
\end{tabular} & -- & \begin{tabular}{c} HT : Superdiffusive\\ LT : Normal \end{tabular} \\
\hline
\end{tabular}
\caption{Diffusivity of the FCLT results in this paper when service times are heavy-tailed (HT, $0 < \alpha < 1$) and light-tailed (LT, $\alpha>1$).}
\label{tab:summary}
\end{table}

The proofs of the FLLNs and FCLTs differ significantly from those of~\cite{PW}, where the approach was to augment the state space by tracking the elapsed service time of the jobs in service. Here, we take a different approach and exploit the fact that conditioned on the stochastic intensity process, the queue can be viewed as an $M_t/G_t/\infty$ queueing system. Thus, under the quenched measure, we use the well known fact~\cite{prekopa,subcritical} that the arrival and service times generate a Poisson random measure on the space $(-\infty,\infty) \times [0, \infty)$ and the state of the queue counts the number of points in rectangles in the space. 
Following classical approach, we then prove the 
FLLN and FCLT by establishing the convergence of the finite dimensional distributions and tightness of the rescaled state process under the quenched regime. 
The proof of tightness is inspired by that of~\cite{RR}, however the bounds now are sample path dependent, which requires some analytical development.

The analysis under the annealed measure is {different in nature}. The proof follows the usual weak convergence recipe~\cite{Bill}. Note that this result, however, does {\it not} follow from~\cite{RR} since the results there crucially depend on the presumption that the re-scaled process is a Poisson random measure; however, this is clearly not the case under the annealed measure, since the performance measures must be averaged over the random environment process. The limits under the subcritical and supercritical regimes are determined by 
differing dominant terms in the re-scaled number-in-system process (in both the quenched and annealed setting). 
We point out in advance here that there are
two stochasticities driving our system, one being the randomness in the
arrival intensity and the other being the service process. Depending on 
the parameter regimes, one can dominate the other.


%
%

A natural question to ask at this juncture is whether the quenched regime is significant at all for performance analysis. As~\cite{ArMoRa2005} observes in the context of Markov models in random environments, {\it "A reasonable first answer to this question could simply be that the quenched approach solves the true question but the averaged \emph{[annealed]} approach, being often much simpler, has the merit of being the first possibility to understand a hard problem..."} Performance analysis is typically an {\it ex ante} exercise, using a stochastic model that predicts the future evolution of the queue. Roughly speaking, an annealed approach to performance analysis permits the prediction of the behavior of the queue  in the {\it typical} random environment. 
The quenched analysis, on the other hand, 
makes predictions for {\it any (almost sure)} 
realization of the random environment. In this sense, our FCLT results 
which show that the limiting $\text{Cox}/G_t/\infty$ queue is closely approximated by an $M_t/G_t/\infty$ queue,  are `robust' to the specific realization of the random environment. As Ben-Arous et al observe~\cite{ArMoRa2005}, this is in some sense the ``true question'' in that, {\it ex ante}, it is unclear what sample path of the random environment will be encountered by the queue, and any analysis should be agnostic to this. 

\subsection{Literature Survey} 
We provide a brief survey of literature relevant to this paper, placing it in context. {Indeed, our paper can be related to multiple threads of research.} First, our main theorems rely on the establishment of a stochastic averaging/homogenization result for the state of the queue, both conditioned and averaged over the underlying stochastic environment model. 
Homogenization/stochastic averaging, has been studied extensively in the 
literature see~\cite{fibich2012averaging,ZHG1,khasminskii2004averaging,khasminskij1968principle,kurtz1992averaging} 
for a small sample of this literature. Typically 
results on stochastic averaging principles in the form of law of large number 
(LLN) involve two processes, one evolving on a fast time-scale and reaching an equilibrium quickly and another one on a slower time-scale that only experiences the former in equilibrium. 
In the limit, the slow process is (typically) approximated by a process with 
``averaged'' coefficients. Averaging principles have been used in finance~\cite{fouque1999financial,fouque2018optimal,blanchet2013continuous}, as well as in heavy-traffic analyses of certain controlled stochastic network models~\cite{coffman1995polling,hunt1994large,perry2011ode}. The averaging principle in our setting is established under the presumption that the process on the fast time-scale is an external stochastic environment, while in~\cite{coffman1995polling,hunt1994large,perry2011ode} the process on the fast time-scale is usually the state of one of the queues in the network, that influences the state of other slowly varying queues. More closely related is the work in~\cite{ZHG1} where the rapidly fluctuating process is an external stochastic environment process.


We also allow for time-dependent service and traffic intensity in our model, and therefore the literature on the analysis of time-varying queues and stochastic models is highly relevant. In the context of the performance analysis of deterministically time-varying queues (not necessarily in a random environment) there is a significant body of work developing both uniform acceleration 
\cite{MMR,honnappa,MM,SS} and many-server heavy-traffic limit theorems 
\cite{lu1,lu2,CH1} to time-varying queues. In much of this literature, the limit processes are shown to be (reflected) fluid or diffusion limit processes. Note that all of this work assumes that the nonstationarity manifests as a deterministic temporal variation. There is also a growing body of work developing asymptotic expansions \cite{ZHG1,ZHG2,MW2,pender1,pender2,KYZ1,KYZ2} of performance metrics. It is well known, however, that traffic arriving at call centers and hospitals displays significant over-dispersion relative to a Poisson process with deterministic intensity~\cite{KW1}, implying that a DSPP is an appropriate model of the traffic in these systems. Much of this literature assumes that either the traffic and/or service processes are Markov modulated, where the underlying stochastic environment process is a finite state Markov chain; the vast majority of the related literature focuses on characterizing stationary behavior, but~\cite{CMRW,pender1} exhibit a couple of examples where asymptotic limit theorems and expansions can be established. The plethora of methods for analyzing time-varying queues will, of course, be crucial for further analyzing the limit $M_t/G_t/\infty$ queue. However, we do not address this fact explicitly in this paper.

Most relevant to our current setting is the expanding literature on infinite server queues in random environments~\cite{OP,HvLM,HM,boxma2019infinite,jansen2019diffusion,dean2020functional,anderson2016functional,blom2013large,blom2013time,blom2014markov,blom2016functional,fralix2009infinite,hellings2012semi,mandjes2016markov,koops2017networks}. More specifically, there are extensive studies on
infinite server queues with Markov-modulated input, in which the arrival 
intensity of the input is $\lambda_i$ with $i$ being the state of the latent 
environment modeled by a Markov jump process.
Results established in this realm touch upon
steady-state, functional central limit theorems and large deviations 
where the latent Markov process is sped up appropriately resulting in a homogenization effect~\cite{anderson2016functional,blom2013large,blom2013time,blom2014markov,blom2016functional,fralix2009infinite,hellings2012semi,mandjes2016markov,jansen2019diffusion}. In all of these papers, the objective 
is an annealed analysis in the sense described above. The more general setting of input models where the arrival intensity changes continuously as a function of the latent environment process has also been extensively studied~\cite{OP,HvLM,HM,boxma2019infinite,dean2020functional}, including homogenization results~\cite{HvLM,HM} showing that the state of the infinite server queue asymptotically aligns with that of a time-homogeneous $M/G/\infty$ queue, steady-state analysis~\cite{OP,boxma2019infinite} and large deviations analysis~\cite{de2014large,jansen2016large,dean2020functional}. \cite{koops2017networks} studies a so-called $M_S/G/\infty$ queue with a Cox arrival process with shot noise intensity and independent and identically distributed (i.i.d.) service times. This intensity model is a special case of the general setting we consider.  This paper establishes both exact and asymptotic results; in the former case, the paper derives an expression for the joint moment generating function of the number in system and arrival intensity process, while in the latter a functional central limit theorem (FCLT) is established in the limit of a large arrival intensity scale. Our annealed results (and, obviously, the quenched results) are fundamentally different since our analysis is one of stochastic homogenization, as opposed to a large scale asymptotic. Note that all of these cited results are in the annealed regime, and to the best of our knowledge the quenched regime has not been studied before for infinite server queues.

\smallskip

\paragraph{\bf Notations}
We introduce a few notation to be used in this paper.
\begin{enumerate}
\item $X \lesssim Y$ means that there exists a constant $0< c < \infty$ such that $X \leq c Y$.
\item $X \gtrsim Y$ means that there exists a constant $0< c < \infty$ such that $X \geq c Y$.
\item $X \asymp Y$ means that there exist two constants $0<c_1 < c_2< \infty$ such that $$c_1 |Y| \leq |X| \leq c_2 |Y| \,.$$
\end{enumerate}

\section{Settings and Results}

As in \cite{HLTY} we model a $\text{Cox}/G_t/\infty$ queue whose arrival 
intensity is thought of as a random function of time 
described by a stochastic process $Z$ whose properties will be described in \cref{hyp:ergodicity}. Hence, our overall model formulation
involves two types of randomness:
one underlying the {\em arrival intensity},
and one coming from the {\em service times}.
Both stochastic ingredients
will be defined on the same probability space $(\Omega, \mathcal F, \PP)$. The expectation with respect to $\PP$ is denoted as $\mathbb E(\cdot)$. In this context, the probability $\PP$ is interpreted as the 
{\em annealed} probability of the random environment. 
In contrast, the {\em quenched} probability, $\PP_Z$,
with a corresponding expectation $\E_Z$,
 is defined by conditioning on the process $Z$. The relation between 
these expression is summarized as
\begin{equation}\label{eqs:quenched-annealed-notation}
\PP_Z (\cdot) = \PP (\cdot | Z) \,, \quad \text{and} \quad \E_Z (\cdot) = \E (\cdot | Z) \,.
\end{equation}
With the above notation, if we further denote the probability distribution of 
$Z$ and its expectation by $\Pbf$ and $\mathbf{E}$, then we have
\begin{equation}\label{def:anneal-P}
\PP = \Pbf \otimes \PP_Z
\quad\text{and}\quad
\mathbb{E}(\cdot) = \mathbf{E}\big[\mathbb{E}_Z(\cdot)\big] \,.
\end{equation}

In this work, we will have statements with respect to both the quenched and
annealed probability. The former means that the statement holds true
for (almost) every realization of $Z$. The following types of convergence 
together with the corresponding terminology and notation will be considered.
\begin{enumerate}
\item	Quenched (vs annealed) almost sure convergence:
{\em $\Pbf$-almost surely}, or equivalently, 
{\em for almost every realization of $Z$}
(vs {\em $\PP$-almost surely}).
\item	Quenched (vs annealed) convergence in probability:
{\em convergence in $\PP_Z$-probability} 
(vs in {\em $\PP$-probability}).
\item	Quenched (vs annealed) convergence in distribution (or
weak convergence):
\begin{equation}\label{def:notation-weak-convg}
\Longrightarrow_{\PP_Z}
\quad\text{vs}\quad
\Longrightarrow_{\PP}.
\end{equation}
\end{enumerate}
We also note that if the underlying randonmess depends only on $Z$, then
$\PP$ and $\mathbb{E}$ are the same as $\Pbf$ and $\mathbf{E}$.

Next we will descibe the assumptions for the randomness used in this paper.

\subsection{Assumptions}

We start by specifying our model for a queueing system in a fast oscillatory random environment. It is based on an ergodic process $Z$ that is incorporated into the arrival rate. The standing assumptions are summarized in the sequel.

\begin{hypothesis}[Ergodicity of Random Environment]\label{hyp:ergodicity}
Let $Z = \{Z_t \, ; \, t \geq 0\}$ be a $\R^d$-valued stochastic process with sample paths that are right continuous with left limits (RCLL) defined on the probability space $(\Omega, \mathcal{F}, \PP)$. The initial distribution $\mathcal{L} (Z_0)$ of $Z$ is denoted by $\rho_0$.
We suppose that $Z$ possesses a unique invariant probability measure $\pi$. 
Furthermore, $Z$ is assumed to be strongly ergodic with 
rate $\frac{1}{2}$ in the following sense: for any $0 < \kappa < \frac{1}{2}$
and any regular enough function $\psi : \R^d \to \R$, there exists a finite 
random variable $C = C_{\kappa, \psi} (Z) > 0$
such that 
for $\Pbf$-almost surely we have
\begin{equation}\label{eq:hyp-ergodic-on-Z}
\left|\frac{1}{t} \iot \psi (Z_u) \, \dd u - \bar{\psi} \right|
\leq \frac{C}{(1+t)^{\kappa}},
\end{equation}
where 
\begin{equation}\label{eq:bar-psi}
\bar{\psi} := \int_{\R^d} \psi(z)\,\pi(\dd z) \,.
\end{equation}
(In the statements and proofs throughout this paper, 
we actually envision $\kappa$ to be as large as possible.
More precisely, we take $\ka = \frac{1}{2} - \delta$, 
with $\delta$ as small as we wish. 
Hence a ``more accurate terminology'' above should be 
``{\em $Z$ is strongly ergodic with rate $\frac{1}{2}^-$}''.)
\end{hypothesis}

\begin{remark}\label{Rmk:Hyp-for-MC-refer}
\cref{hyp:ergodicity} is satisfied for a wide class of reversible Markov chains, as well as for Langevin type stochastic differential equations driven by Gaussian processes. We refer to \cite{HLTY} Remarks 2.2 and 2.3 for further details. On the other hand, our constraint $\ka < 1/2$ is due to the fact that one cannot reach $\ka \geq 1/2$ in our main application contexts, as dictated by the Law of Iterated Logarithm. 
\end{remark}

With the above, we now describe the arrival process which incorporates
a random intensity and scaling to a heavy traffic regime. 

\begin{hypothesis}[Arrival Intensity]\label{hyp:arr-times-Gamma}
Given an oscillating parameter $\ep > 0$ and a scaling exponent $\beta \geq 0$, the sequence of arrival times $\{\Gamma_k^{\ep} \, ; k \geq 1 \}$ is distributed as a Poisson process with rescaled nonhomogeneous intensity 
\begin{equation}\label{eq:mu-ep-beta}
\mu^{\ep} (\dd s) := \frac{1}{\ep^{\beta}} \la (s) \psi (Z_{s/\ep}) \dd s, \quad \text{for} \,\, s \geq 0 \,.
\end{equation}
where $Z$ fulfills \cref{hyp:ergodicity} and
$\lambda : [0,\infty) \to \R_+$, $\psi : \R^d \to \R_+$ are positive and bounded Lipschitz (or $C^1$-) functions.
\end{hypothesis}

\begin{remark}
Notice that expression \eqref{eq:mu-ep-beta} incorporates two phenomena.
First, the factor $\ep^\beta$ accounts for a heavy traffic parameter in our model. 
Second, the $\ep$-scaling in $Z_{s/\ep}$ is our fast oscillating feature for the arrival rate. It potentially leads to homogenization features.
\end{remark}

Next we turn to the service times for our queueing system.
These are given by a family $\{L_k^{\ep} \,;\, k \geq 1\}$ of random variables which are independent conditional on the arrival process $\{\Gamma_k^{\ep} \, ; k \geq 1 \}$. Their laws are denoted by 
$\mathcal{L} (L_k^{\ep} | \Gamma_k^{\ep}) \equiv \nu (\Gamma_k^{\ep}, \dd r)$ for all $k$. Here we consider
\begin{equation}\label{set:nu-s-dr}
\nu=\left\{\nu(s,\dd r) \,;\, s \geq 0 \,, r \geq 0\right\}
\end{equation}
as a family of conditional regular laws. With these, the complementary cumulative distribution function of the service times $\bar{F}_s (r)$ hereby can be written as
\begin{equation}\label{eq:def-bar-F}
\bar F_s (r) := \int_r^\infty \nu(s,\dd \tau),
\quad r \geq 0 \,.
\end{equation}
The main assumption on the measures $\nu$ is given below.

\begin{hypothesis}[Service time distribution]\label{hyp:F-bar-increments}
Let $\nu$ be the family of measures defined by \eqref{set:nu-s-dr}. We suppose that every $\nu (s , \cdot)$ admits a density $\ell_s$, that is, $\nu (s, \dd r) = \ell_s (r) \dd r$. Moreover,
there exist $\al > 0$ and constants $c,C > 0$ such that the 
family $\ell_s$ verifies
	\begin{eqnarray}\label{def:F-bar-generalization}
		c  \left( \frac{1}{r^{1+\al}} \wedge 1\right)\,\, \leq \,\,\,\,\,\ell_s (r) &\leq& C \left( \frac{1}{r^{1+\al}} \wedge 1\right) ,\quad \text{for all} \quad r, s >0 \,;\\
\label{ineq:partial-ell-s-bounds}
\frac{\partial \ell_s (r)}{\partial s} &\leq& C \left( \frac{1}{r^{1+\al}} \wedge 1\right)\, , \quad \text{for all} \quad r, s >0 \,.
	\end{eqnarray}
\end{hypothesis}

\begin{remark}
Note that for $0 < \alpha < 1$, we have $EL_k^\ep = \infty$ 
while for $\alpha > 1$, we have $EL_k^\ep < \infty$. Hence we call the former 
regime `heavy-tailed' while the latter `light-tailed'. If $\alpha > 2$, 
then we would further have $\text{Var}(L_k^\ep) < \infty$.
\end{remark}

With the above notation and hypotheses in hand, our queueing system is classically described by a point process. Namely for $\ep >0$, we define the following counting measure on $\R_+ \times \R_+$:
\begin{equation}\label{eq:M-sum1}
M^{\ep} := \sum_{k=1}^{\infty} \delta_{(\Gamma_k^{\ep}, L_k^{\ep})} \,.
\end{equation}
Our main quantity of interest is the number of active jobs $N^\ep(t)$ at time 
$t$ which is derived from $M^{\ep}$ by the following formula:
\begin{equation}\label{eq:N-t1}
N^{\ep} (t) = \sum_{k=1}^{\infty} \mathbbm{1}_{\{\Gamma_k^{\ep} < t< \Gamma_k^{\ep} + L_k^{\ep} \}} = M^{\ep} \Big(\{(x,y) \in \mathbb{R}_+ \times \mathbb{R}_+\,,\, x< t< x+y\} \Big) \,.
\end{equation}
We emphasize here that the quantities $\mu^\ep$, $M^\ep$ and $N^\ep$ all
depend on a realization of $Z$. Quenched statements refer to the case that
these random variables are conditioned on a given fixed $Z$.

For the purpose of analysis in the annealed 
setting, we 
introduce the following $L^1$-type ergodicity condition.

\begin{hypothesis}[$L^1$-ergodicity of random environment]
\label{hyp:natual-annealed-bound}
Let $Z$ be the $\R^d$-valued process and $\psi$ the function introduced in \cref{hyp:ergodicity}. Also recall the $\bar\psi$ defined by \eqref{eq:bar-psi}. Then we assume that 
\begin{equation}\label{ineq:hyp-annealed-psi}
\lim_{t\to\infty}
\mathbb{E}\left|
\frac{1}{t} \iot \psi (Z_r) \, \dd r - \bar{\psi} \right|
=
0.
\end{equation}
\end{hypothesis}
Observe that the above condition will certainly follow if we assume
some uniform integrability of the constant $C$ (which in fact is a
random variable depending on $Z$) in \cref{hyp:ergodicity}.
As noted in the introduction, since the annealed analysis demonstrates the limiting behavior of the infinite server queue in a ``typical'' stochastic environment, hence it is natural that a somewhat
weaker ergodicity condition can be assumed.

Furthermore, for the analysis of CLT type fluctuation phenomena, we will make
the following hypothesis about the process $Z$.

\begin{hypothesis}[FCLT for $Z$]\label{hyp:clt-Z}
Let $Z$ be the RCLL stochastic process defined in \cref{hyp:ergodicity} and $\psi : \R^d \to \R_+$ a positive Lipschitz function as in \cref{hyp:arr-times-Gamma}. In addition, we consider $\psi \in L^2 (\R)$. 
For $\ep >0$, define a process $Y^{\ep} = \{Y^{\ep} (t) \,;\, t \geq 0\}$ by
\begin{equation*}
Y^{\ep} (t) = \frac{1}{\sqrt{\ep}} \iot \Big(\psi (Z_{s/\ep}) - \bar\psi \Big) \,\dd s \,.
\end{equation*}
Then we assume that there exists a constant $\sigma_{\psi} >0$ and a Brownian motion $W$ defined on $(\Omega, \mathcal{F}, \PP)$ such that the following convergence in law holds true in the space $C(\R_+)$ of continuous functions (equipped with the topology of uniform convergence on compact sets):
\begin{equation}\label{lim:Y-ep-convg-BM}
\{Y^{\ep} (t) \,;\, t \geq 0\} \Longrightarrow_{\PP} \{\sigma_{\psi} W (t)\,;\, t \geq 0\} \,.
\end{equation}
\end{hypothesis}


\begin{remark}
Central limit theorems like \eqref{lim:Y-ep-convg-BM} for Markov processes are treated at length in \cite{KLO}, subsequent to the celebrated paper \cite{KV}. A typical example can be found in \cite[Theorem 2.7]{KLO}. We also refer to \cite{CCG} for Multi-times Central Limit Theorems (MCLT) for Markov processes. In a non-Markov setting, let us quote the continuous functional Breuer-Major theorem spelled out in \cite{CaNua}. This result is valid for stationary Gaussian processes with a mild integrability assumption on its covariance function. Eventually, observe that \cref{hyp:clt-Z} is used in other queueing theory references, see e.g. \cite{PW}.
\end{remark}

\subsection{Different regimes for the rescaled Queue}\label{subsect:intro-regimes-rescaledqueue}

Our previous article \cite{HLTY} was conceived with a $\text{Cox}/G_t/\infty$ queue with a random environment whose intensity is given by \eqref{eq:mu-ep-beta}, albeit in a situation for which $\beta = 0$. 
For convenience, 
we use $\mu_0^\ep, M^\ep_0$ and $N^\ep_0$ to denote
the corresponding quantities defined above, specific to the case 
$\beta=0$. More precisely,
\begin{enumerate}
\item
$\displaystyle
\mu^\ep_0(\dd s) := \lambda(s)\psi(Z_{s/\ep}) \dd s
$ is the arrival intensity;

\item
$M^\ep_0$ is a Poisson random measure 
with mean measure:
$$
\nu^\ep_0(\dd x, \dd y) = \nu(x,\dd y)\mu^\ep_0(\dd x), 
$$

\item	
$N^\ep_0(t)$ is a a Poisson random variable with parameter:
\begin{equation}\label{eq:m-ep-0-defined-int}
m^\ep_0(t) = \int_{(x,y): x< t< x+y}\nu(x,\dd y)\mu^\ep_0(\dd y)
= \int_0^t\int_{t-x}^\infty\nu(x,\dd y)\mu^\ep_0(\dd y).
\end{equation}
\end{enumerate}

The main results in \cite{HLTY} give a description of the limiting behavior of the random quantities
$\mu^\ep_0, M^\ep_0, N^\ep_0$. Namely, the measure $\mu^\ep_0$ converges to $\bar{\mu}_0$, where 
$$
\bar{\mu}_0 (\dd s) = \bar\psi \la (s) \dd s \,.
$$
Next, the conditional Poisson random measure $M^\ep_0$ also converges to 
$\bar{M}_0$ which is a Poisson measure with intensity
\begin{equation}
\bar{\nu}_0(\dd x, \dd y)
=\bar{\psi} \,\bar{F}_x(y) \,\la (y)\, \dd x \,\dd y \,. \label{nu-0-def}
\end{equation}
Eventually $N^\ep_0$ converges to a process with Poisson marginals and mean process
\begin{equation}
\bar{m}_0(t) =\E_Z [\bar N_0 (t)] = \Lambda (t)\bar{\psi} \,,
\label{m0-def}
\end{equation}
where
\begin{equation*}
\Lambda (t) =\int_0^t h (s,t) \, \dd s ,
\quad \text{with} \quad h(s,t) = \la (s) \bar{F}_s(t-s) \,,
\end{equation*}
with $\bar\psi$ and $\bar F$ given by \eqref{eq:bar-psi} and
\eqref{eq:def-bar-F}, respectively.

In \cite{HLTY}, the convergence of $N^\ep_0$ to $\bar{N}_0$ is established in 
the sense of stochastic process in both quenched and annealed regimes.
We notice that the limiting process $\bar{N}_0$ is still a Poisson point process, which is in sharp contrast with the rescaled version presented in the current contribution.
More specifically and as mentioned in the introduction, the present work considers the case $\beta>0$
corresponding to a 
sped-up arrival process leading to a heavily increased traffic regime.
Hence, we need to rescale the number of active jobs $N^\ep(\cdot)$ in order
to obtain a meaningful limit. 
We first  heuristically describe our results with precise statements 
deferred to \cref{Sect:LLN-for-N-ep,sect:clt-finite-dim-dist,Sect:annealed-analysis}.

Note that under \cref{hyp:arr-times-Gamma,hyp:F-bar-increments}, 
when conditioned on $Z$, 
the random variable $N^\ep(t)$ defined by \eqref{eq:N-t1}
is again a Poisson random variable (conditionally on $Z$) with parameter 
$m^\ep(t) = \frac{1}{\ep^\beta}m^\ep_0(t)$.
Hence, the LLN for $N^\ep(t)$ can be stated as:
\begin{equation}\label{eqlim:lln-N-ep-over-m-ep}
\frac{N^\ep(t)}{m^\ep(t)}=
\frac{\ep^\beta N^\ep(t)}{m^\ep_0(t)}
\longrightarrow 1\,.
\end{equation}
We notice that the above assertion is valid under the quenched probability $\PP_Z$.

The next step is to analyze the fluctuation phenomena of $N^\ep (t)$.
Since $m^{\ep} (t) = \frac{m^\ep_0(t)}{\ep^\beta}\rightarrow\infty$, 
from the elementary property of Poisson random variables, we have that,
\begin{equation}\label{a1}
\frac{N^\ep(t) - \frac{m^\ep_0(t)}{\ep^\beta}}{\sqrt{\frac{m^\ep_0(t)}{\ep^\beta}}}\Longrightarrow_{\PP_Z} \cn(0,1) \,,
\end{equation}
where we recall our convention \eqref{def:notation-weak-convg} on quenched convergence in distribution.
However, it is more natural to consider the following rescaled version which is 
centered around a {\em deterministic quantity}:
\begin{equation}\label{centered-Nep}
G^{\ep} (t) :=
\ep^{\frac{\beta}{2}}\Big(N^\ep(t) - \frac{\bar{m}_0(t)}{\ep^\beta}\Big).
\end{equation}
The above expression can be decomposed into
{\em two types of fluctuations}:
\begin{equation}\label{clt-decompose}
G^{\ep} (t) := G^{\ep}_1 (t) + G^{\ep}_2 (t)
\end{equation}
where 
\begin{equation}\label{eqdef:Q-ep-1,2-fluctuation}
G^{\ep}_1 (t) =
\frac{\ep^\beta N^\ep(t) - m^\ep_0(t)}
{\ep^\frac{\beta}{2}}
\quad \text{and} \quad
G^{\ep}_2 (t) = \ep^{\frac{1- \beta}{2}}\Big(\frac{m^\ep_0(t)
-\bar{m}_0(t)}{\sqrt{\ep}}\Big) \,.
\end{equation}
Observe that $G^{\ep}_1 (t)$ is dictated by the fluctuation of the service times while $G^{\ep}_2 (t)$ comes from the driving force for the arrival process. The limiting behavior in~\eqref{clt-decompose}
is thus determined by the {\em dominating term}.
In fact, under \cref{hyp:clt-Z}, we have heuristically,
\begin{equation*}
\frac{m^\ep_0(t)-\bar{m}_0(t)}{\ep^\frac{1}{2}}
\Longrightarrow_{\PP}
\cn(0, \bar\sigma^2 (t)) \,,
\end{equation*}
where the expression for $\bar\sigma^2 (t)$ will be established in \cref{thm-FCLT-super-univeriate-annealed}.
Our main message is that for $\beta < 1$, the limit for \eqref{clt-decompose} is described by the asymptotic behavior of $G_1^{\ep} (t)$, namely 
\begin{equation}
G^{\ep} (t) =
\ep^{\frac{\beta}{2}}\Big(N^\ep(t) - \frac{\bar{m}_0(t)}{\ep^\beta}\Big)
\sim 
G^{\ep}_1 (t) =
\frac{\ep^\beta N^\ep(t) - m^\ep_0(t)}
{\ep^\frac{\beta}{2}} \sim O (1) \,,
\label{sub-heuristic}
\end{equation}
while for $\beta > 1$, the limit is described by $G^{\ep}_2 (t)$. In the latter case,
we need to consider the following rescaling in order to obtain a
meaningful limit:
\begin{equation}
\ep^{\beta-\frac{1}{2}}
\Big(N^\ep(t) - \frac{\bar{m}_0(t)}{\ep^\beta}\Big)
\sim
\frac{m^\ep_0(t)-\bar{m}_0(t)}{\sqrt{\ep}} \sim O(1).
\label{super-heuristic}
\end{equation}

Based on the above, we make the following definition.
\begin{definition}
We call $\beta < 1$, $=1$, and
$>1$ as {\em subcritical, critical}, and 
{\em supercritical} regimes. 
\end{definition}

Some remarks are in place. 
It is apparent that in the supercritical regime, a quenched FCLT does not exist, as a meaningful limit for $G^\ep_2$ cannot be established under the quenched measure. 
For then, one would have to argue for the existence of a 
{\it pathwise} limit under the CLT scaling for $Z$.
 Note that we have not treated the critical case in this article for the sake of conciseness. However, our conjecture is that it would yield an annealed CLT where $G^\ep_1$ and $G^\ep_2$ would give contributions of equal magnitude. 
We also recall the discussion in the introduction section 
that this scenario in fact
parallels the conclusions in \cite[Section 2.2 and 2.3]{Peterson08}
for random walks in random environments.

\subsection{Outline of Main Results}
The main contribution in this paper is the development of  (functional) central limit theorems 
in both the quenched and annealed setting for a rescaled version of the number 
of active jobs $N^{\ep}$ defined by~\eqref{eq:N-t1}. As explained above,
the precise statement and formulation will depend on the value of $\beta$. Here is a road map of the main results established :

\begin{enumerate}[leftmargin=*]
\item Quenched Analysis (Section~\ref{Sect:quenched-analysis})
\begin{enumerate}
\item	Quenched LLN (sub- and supercritical) 
Theorem \ref{thm-LLN-quenched} 
\item	Quenched CLT (subcritical case)
Theorem \ref{thm-FCLT-quenched}
\end{enumerate}
\item Annealed Analysis (Section~\ref{Sect:annealed-analysis})
\begin{enumerate}
\item	Annealed LLN (sub- and supercritical)
Theorem \ref{thm-LLN-annealed}
\item	Annealed CLT (subcritical case)
Theorem \ref{thm-FCLT-sub-annealed}
\item	Annealed CLT (supercritical case)
Theorem \ref{thm-FCLT-super-annealed}
\end{enumerate}
\end{enumerate}

We close this section with a few words about our techniques of proof. The basic ingredient for the finite dimensional convergence in distribution of $G^{\ep}$ is based on standard considerations for convergence of Poisson random variables. However, this has to be combined with a thorough and delicate inspection of the almost sure behavior of our random parameters. This almost sure behavior relies on the elements of ergodic theory. It should also be noticed that the functional part of our CLT hinges on intricate arguments about modulus of continuity for $G^{\ep}$. 


\section{Quenched Analysis}~\label{Sect:quenched-analysis}

\subsection{Quenched Law of Large Numbers for $N^\ep(t)$}\label{Sect:LLN-for-N-ep}
In this section, we analyze the first-order asymptotic behavior of the process $N^{\ep}$ defined by \eqref{eq:N-t1}. 
The statements to be established are
applicable for both the subcritical and supercritical regimes.
We emphasize here again the notion that ``quenched’’ means that the 
statements depend on a given, fixed realization of $Z$.

We start by recalling some
established facts about the mean measure of the Poisson point process
$M^{\ep}$ and their limiting description as $\ep\longrightarrow0$. Specifically, we first give an expression for the mean value of $N^{\ep} (t)$. 
The proof is omitted as it is based on classical considerations on Laplace transforms for Poisson measures; see e.g. \cite{Cin}.

\begin{proposition}\label{prop:m-t-of-N-t}
Let $M^\ep$ be defined by \eqref{eq:M-sum1} and $\beta \geq 0$.
Then under the quenched probability 
$\mathbb{P}_Z$,
$M^\ep$ is a Poisson random measure with mean measure given by
\begin{equation}\label{eq:def-mu-la-psi}
{\nu}^{\ep} (\dd x, \dd y) = \nu (x, \dd y) \mu^{\ep} (\dd x) \,,
\quad \text{with} \quad
\mu^{\ep} (\dd s) = \frac{1}{\ep^{\beta}} \la (s) \psi (Z_{s/\ep}) \dd s \,,
\end{equation}
where $\nu$ is introduced in \eqref{set:nu-s-dr}.
Furthermore,  we have that for any $t > 0$, 
$N^{\ep} (t)$ is a Poisson random variable with parameter $m^{\ep} (t) = \E_Z [N^{\ep} (t)]$ given by
\begin{equation}\label{eq:m-t-def}
m^{\ep} (t) = \int_{\{(x,y): x<t< x+y\}} \nu (x, \dd y) \mu^{\ep} (\dd x) = \iot \int_{t-x}^{\infty} \nu (x, \dd y) \mu^{\ep} (\dd x).
\end{equation}
\end{proposition}

\noindent
We give an alternative representation for further use. It is readily checked from \eqref{eq:m-t-def} that $m^{\ep} (t)$ can be recast as
\begin{equation}\label{def:m-ep-expansionwith-m-ep-0}
m^{\ep} (t) = \frac{1}{\ep^\beta} \iot h (s,t)\, \psi (Z_{s/\ep}) \,\dd s = \frac{1}{\ep^\beta}\, m^{\ep}_0 (t) \,,
\end{equation}
where $m^{\ep}_0 (t)$ is defined by \eqref{eq:m-ep-0-defined-int} and where recalling from \eqref{m0-def}, we have set $h(s,t) = \la (s) \bar F_s (t-s)$.

Next we establish the limiting behavior of $m^{\ep} (t)$.

\begin{proposition}
\label{prop:gen-limit-m-ep}
Let us assume the same set-up and notations as in 
\cref{prop:m-t-of-N-t}. Furthermore,
suppose that \cref{hyp:ergodicity,hyp:arr-times-Gamma,hyp:F-bar-increments} are satisfied. 
Then the following two statements hold $\Pbf$-almost surely. 

\begin{enumerate}[wide, labelwidth=!, labelindent=0pt, label=\emph{(\roman*)}]
\setlength\itemsep{.1in}

\item\label{it:convergence-nep-i}
Let $\beta = 0$. For any $t> 0$, 
\begin{equation}\label{eq:lim-m-ep}
\lim_{\ep \to 0}  m^{\ep}_0 (t) = \bar{m}_0 (t) := \Lambda (t) \bar{\psi} \, .
\end{equation}
(For convenience, we repeat here the definitions 
\eqref{eq:bar-psi}, \eqref{m0-def} of the following quantities:
\begin{equation}\label{eq:quantity-sigma-psi-bar}
\bar\psi = \int_{\R^d} \psi(z)\,\pi(\dd z) 
\quad \text{and} \quad
\Lambda (t) = \iot \la (s) \bar{F}_s (t-s) \, \dd s = \iot h (s,t) \dd s \,,
\end{equation}
where $\pi$ is the invariant measure of the process $Z$ introduced in 
\cref{hyp:ergodicity}.)

\item\label{it:convergence-nep-ii}
Let $\beta > 0$. For any $t>0$,
\begin{equation}\label{eq:limit-epbeta-m-ep}
\lim_{\ep \to 0} \ep^{\beta} m^{\ep} (t) =\bar{m}_0 (t) 
= \Lambda (t) \bar\psi \,.
\end{equation}
In particular, $\Pbf$-almost surely,
\begin{equation}\label{eq:limit-rescaled-m-ep}
\lim_{\ep \to 0} m^{\ep} (t) = \infty \,.
\end{equation}
\end{enumerate}
\end{proposition}

\begin{proof}
With the expression \eqref{def:m-ep-expansionwith-m-ep-0} for $m^{\ep} (t)$, 
our claims \eqref{eq:lim-m-ep} and \eqref{eq:limit-epbeta-m-ep} are easy consequences of the general \cref{weighted.ergodic} below.
\end{proof}

We now present a useful averaging lemma which is crucial to complete
the proof of \cref{prop:gen-limit-m-ep}. It will also be used in several other steps throughout the paper.

\begin{lemma}[Weighted ergodicity]\label{weighted.ergodic}
Assume \cref{hyp:ergodicity} holds true and let $Z$ be the process involved in this assumption. Consider a $C^1$ and bounded function $h$ defined on $[0,T]^2$.
Then the following statements hold for $\Pbf$-almost surely.
\begin{equation}\label{homog.limit}
\lim_{\ep\to0}\int_0^t h(s,t)\,\psi(Z_{s/\ep})\, \dd s
=\bar{\psi}\int_0^t h(s,t)\,\dd s \,.
\end{equation}
Moreover, if $\beta < 2\kappa < 1$, i.e. {\em sub-critical case}, then
\begin{equation}\label{homog.limit.subrate}
\lim_{\ep\to0}
\frac{1}{\ep^{\beta/2}}
\left|\int_0^t h(s,t)\,\psi (Z_{s/\ep})\,\dd s
-\bar{\psi}\int_0^t h(s,t)\,\dd s\right| = 0 \,.
\end{equation}
In particular, by setting
$\displaystyle h(s,t) = \la (s) \bar F_s (t-s)$, we have
\begin{equation}\label{m.ep.rate.subcrit}
\lim_{\ep \to 0} m^{\ep}_0 (t) = \bar m_0 (t) \,
\quad\text{and}\quad
\lim_{\ep \to 0}
\frac{1}{\ep^{\beta/2}}
\left| m^{\ep}_0 (t) - \bar m_0 (t)\right|
= 0 \,.
\end{equation}
\end{lemma}

\begin{proof}
We first prove \eqref{homog.limit}. A simple integration by parts procedure shows that
\begin{eqnarray}
&&\int_0^t h(s,t) \,\psi (Z_{s/\ep})\,\dd s \notag\\
& = & 
\int_0^t h(s,t) \, \frac{\dd}{\dd s} \left(\int_0^s\psi(Z_{r/\ep})\,\dd r\right)\notag \\
& = & 
h(t,t)\int_0^t\psi(Z_{r/\ep})\,\dd r
-
\int_0^t \frac{\dd}{\dd s}h(s,t)
\left(\int_0^s\psi(Z_{r/\ep})\,\dd r\right)\,\dd s \notag \\
& = & 
h(t,t)t\left[\frac{1}{t}\int_0^t\psi(Z_{r/\ep})\,\dd r
-\bar{\psi}\right] + th(t,t)\bar{\psi}\notag \\
& & 
-
\int_0^t \frac{\dd}{\dd s}h(s,t) \,s\left[
\frac{1}{s}\int_0^s\psi(Z_{r/\ep})\,\dd r - \bar{\psi}\right]
\,\dd s
- \bar{\psi} \,\iot \frac{\dd}{\dd s}h(s,t) \,s\,\dd s \,. \label{eq:compute-hpsi-Z-int}
\end{eqnarray}
Moreover, another integration by parts reveals that
\begin{eqnarray*}
th(t,t)\bar{\psi} - \bar{\psi}\, \iot \frac{\dd}{\dd s}h(s,t) \,s\,\dd s
= \bar{\psi}\int_0^t h(s,t)\, \dd s \,.
\end{eqnarray*}
Plugging this information into \eqref{eq:compute-hpsi-Z-int}, we obtain
\begin{eqnarray}
&&\int_0^t h(s,t) \,\psi (Z_{s/\ep})\,\dd s - \bar{\psi}\int_0^t h(s,t)\, \dd s \notag\\
&=& h(t,t)t\left[\frac{1}{t}\int_0^t\psi(Z_{r/\ep})\,\dd r
-\bar{\psi}\right] - \int_0^t \frac{\dd}{\dd s}h(s,t) \, s\left[
\frac{1}{s}\int_0^s\psi(Z_{r/\ep})\,\dd r - \bar{\psi}\right]
\,\dd s \,. \label{eq:compute-hpsi-Z-int-var}
\end{eqnarray}
Now according to \cref{hyp:ergodicity}, we have (in the $\Pbf$-almost sure sense) that for a finite random variable $C = C_{\kappa, \psi} (Z) > 0$,
\begin{eqnarray*}
\left|
\frac{1}{s}\int_0^s\psi(Z_{r/\ep})\,\dd r - \bar{\psi}
\right|
\leq
\frac{C}{(1+\frac{s}{\ep})^\ka} \,.
\end{eqnarray*}
Hence for a fixed $t$, we easily get
\begin{equation}\label{ineq:lim-h-psi-Z-compute}
\lim_{\ep\to0}\left|h(t,t)t\left[\frac{1}{t}\int_0^t\psi(Z_{r/\ep})\,\dd r
-\bar{\psi}\right]\right|
\leq\lim_{\ep\to0}\left|h(t,t)t\frac{C}{(1+\frac{t}{\ep})^\ka}
\right|
=
0 \,.
\end{equation}
On the other hand, due to the regularity properties of $h$, one can apply Lebesgue Dominated Convergence Theorem to get
\begin{multline}\label{ineq:LDCT-lim-int-h-psi}
\lim_{\ep\to0}\left|\int_0^t \frac{\dd}{\dd s}h(s,t) s\left[
\frac{1}{s}\int_0^s\psi(Z_{r/\ep})\,\dd r - \bar{\psi}\right]
\,\dd s\right| \\
\leq
\lim_{\ep\to0}
\int_0^t 
\left|\frac{\dd}{\dd s}h(s,t) s\right|
\frac{C}{(1+\frac{s}{\ep})^\kappa}
\,\dd s
=0 \,.
\end{multline}
Gathering \eqref{ineq:lim-h-psi-Z-compute} and \eqref{ineq:LDCT-lim-int-h-psi} into \eqref{eq:compute-hpsi-Z-int-var} and taking limit complete the proof of \eqref{homog.limit}.

We next prove \eqref{homog.limit.subrate}.
Note that owing to the fact that $\beta < 2 \ka$, we have
\begin{eqnarray*}
\lim_{\ep\to0}\frac{1}{\ep^{\beta /2}}
\left|
\frac{1}{s}\int_0^s\psi(Z_{r/\ep})\,\dd r - \bar{\psi}
\right|
\leq \lim_{\ep\to0}\frac{C}{\ep^{\beta /2}(1+\frac{s}{\ep})^\ka}
\leq \lim_{\ep\to0}\frac{C\ep^{\ka - \beta /2}}{s^\ka}
=0 \,.
\end{eqnarray*}
As $\kappa < \frac{1}{2} < 1$, $s^{-\kappa}$ is integrable. We can then
again apply Lebesgue Dominated Convergence Theorem as in \eqref{ineq:LDCT-lim-int-h-psi} to get the result.
Then statement \eqref{m.ep.rate.subcrit} is a direct consequence of~\eqref{homog.limit.subrate}
in the specific case $h(s,t) = \la (s) \bar F_s (t-s)$.
\end{proof}

The following Lemma is useful in stating 
results in the quenched setting which are required to hold true
$\Pbf$-almost surely for all parameters.
\begin{lemma}[$\Pbf$-almost sure limit with parameter]\label{limit.parameter}
Let $\{f_\ep: \Theta\times\Omega\longrightarrow\R\}_{\ep > 0}$ 
where $\Theta$ is an open subset of $\R$. Suppose the following properties 
hold:
\begin{enumerate}
\item	For $\Pbf$-almost every $\omega$, $f_\ep(\cdot, \omega)$ 
is Lipschitz with Lipschitz constant $L(\omega)<\infty$ uniformly in $\ep >0$, i.e.
\begin{equation}
|f_\ep(\theta_1, \omega) - f_\ep(\theta_2,\omega)|
\leq L(\omega)|\theta_1 - \theta_2|
\quad\text{for any $\theta_1,\theta_2\in\Theta$.}
\end{equation}
\item	For each $\theta\in\Theta$, for $\Pbf$-almost every $\omega$,
the limit $\lim_{\ep\to0}f_\ep(\theta,\omega)$ exists and is denoted by 
$g(\theta,\omega)$.
\end{enumerate}
Then there is a version $\tilde{g}:\Theta\times\Omega\longrightarrow\mathbf{R}$
of $g$\footnote{That is, for any $\theta\in\Theta$, 
$\Pbf\{\omega:\tilde{g}(\theta,\omega)=g(\theta,\omega)\}=1$.}
such that 
for $\Pbf$-almost every $\omega$, we have 
\begin{equation}
\lim_{\ep\to0}f_\ep(\theta,\omega) = \tilde{g}(\theta,\omega)
\quad\text{for all $\theta\in\Theta$}.
\end{equation}
(Note that the $\omega$ works {\em for all $\theta$}.)
\end{lemma}
\begin{proof}
Without loss of generality, there is an $\Omega^\circ\subset\Omega$ 
with full $\Pbf$-measure such that for all $\omega\in\Omega^\circ$, we have
$\lim_{\ep\to0}f_\ep(\theta,\omega) = g(\theta,\omega)$ for $\theta$ in 
a countable dense subset $\Theta^\circ$ of $\Theta$ and 
$f_\ep(\cdot, \omega)$ is Lipschitz with constant $L(\omega)$.
Now we define $\tilde{g}:\Theta\times\Omega^\circ\longrightarrow\mathbb{R}$,
\[
\tilde{g}(\theta,\omega) = \left\{
\begin{array}{cll}
g(\theta,\omega) & \text{for} & \theta\in\Theta^\circ;\\
\lim_{\theta_i\in\Theta^\circ\to\theta}g(\theta_i,\omega) & \text{for} &
\theta\in\Theta\backslash\Theta^\circ.
\end{array}
\right.
\]
By the uniform Lipschitz continuity of $f_\ep$ in $\theta$, the function
$\tilde{g}$ is well-defined and also Lipschitz with constant $L(\omega)$. Now
for any $\omega\in\Omega^\circ$ and $\theta\in\Theta\backslash\Theta^\circ$,
we have,
\begin{eqnarray*}
& & |f_\ep(\theta,\omega) - \tilde{g}(\theta,\omega)|\\
& \leq & 
|f_\ep(\theta,\omega) -f_\ep(\theta_i,\omega)| 
+ |f_\ep(\theta_i,\omega) - \tilde{g}(\theta_i,\omega)|
+ |\tilde{g}(\theta_i,\omega)-\tilde{g}(\theta,\omega)|
\end{eqnarray*}
where $\theta_i\in\Theta^\circ\longrightarrow\theta$. Note that all the three
terms in the above can be made as small as desired due to the fact that
$\theta_i\longrightarrow\theta$ and $\ep\longrightarrow0$. Thus, the result follows.
\end{proof}

We now state a law of large numbers for the process $N^{\ep}$ in the rescaled version when our parameter $\beta$ is strictly positive.

\begin{theorem}\label{thm-LLN-quenched}
Let $N^{\ep}$ be the process defined by \eqref{eq:N-t1}, and assume that \cref{hyp:ergodicity,hyp:arr-times-Gamma,hyp:F-bar-increments} are satisfied. Furthermore, suppose $\beta > 0$. Then $\Pbf$-almost surely, we have
for any fixed $t\geq 0$ that the following limit holds as $\ep\to 0$,
\begin{equation}\label{eq:LLN-N-ep-equal-1}
\frac{N^{\ep} (t)}{m^{\ep} (t)} \longrightarrow 1
\quad\text{in $\PP_Z$-probability} \,,
\end{equation}
where $m^{\ep} (t)$ is defined by \eqref{eq:m-t-def}. 
More precisely, $\Pbf$-almost surely we have that for all $\delta > 0$,
\begin{equation*}
\lim_{\ep \to 0} \, \PP_Z \left(\left|\frac{N^{\ep} (t)}{m^{\ep} (t)} - 1 \right| \geq \delta \right) = 0 \,.
\end{equation*}
In addition, as $\ep\to 0$, the following convergence holds true:
\begin{equation}\label{NepLimitm0}
\frac{\ep^{\beta} N^{\ep} (t)}{\bar{m}_0(t)} \longrightarrow 1
\quad\text{in $\PP_Z$-probability} \,.
\end{equation}
\end{theorem}

\begin{proof}
For convenience we set
\begin{equation}\label{eq:U-ep-fraction}
U^{\ep} (t) = \frac{N^{\ep} (t)}{m^{\ep} (t)} \,.
\end{equation}
Then the result \eqref{eq:LLN-N-ep-equal-1} is equivalent to $\lim_{\ep \to 0} U^{\ep} (t) = 1$ in 
$\PP_Z$-probability. Since the limit is a constant, the statement follows from a mere 
convergence in distribution. This can be characterized in terms of 
Laplace transform: we will show that $\Pbf$-almost surely, 
the following holds true for all $\theta \geq 0$:
\begin{equation}\label{eq:limit-exp-theta-U-ep}
\lim_{\ep \to 0} \E_Z \left[e^{-\theta U^{\ep} (t)} \right] = e^{-\theta} \,.
\end{equation}

First, consider a fixed $\theta$. Observe that according to 
\cref{prop:m-t-of-N-t}, $N^{\ep} (t)$ is a Poisson random variable with parameter $m^{\ep} (t)$. Hence, the Laplace transform of $U^{\ep} (t)$ defined in~\eqref{eq:U-ep-fraction} is given by
\begin{equation}\label{Laplacetransform:N-ep/m-ep}
\E_Z \left[e^{-\theta U^{\ep} (t)} \right]
=
\exp \left\{ m^{\ep} (t) (e^{-\theta /m^{\ep} (t)} - 1) \right\} .
\end{equation}
In order to take limits in \eqref{Laplacetransform:N-ep/m-ep}, we first compute the limit for the exponent on the right hand side. Thanks to the facts 
$e^{-u} - 1 = -u + o (u)$ as $u \to 0$ and 
$\lim_{\ep \to 0} m^{\ep} (t) = \infty$, we get
\begin{equation*}
\lim_{\ep \to 0} m^{\ep} (t) (e^{-\theta /m^{\ep} (t)} - 1)
=
\lim_{\ep \to 0} m^{\ep} (t) \, \frac{-\theta}{m^{\ep} (t)}
= - \theta .
\end{equation*}
Plugging this relation into \eqref{Laplacetransform:N-ep/m-ep}, we thus have,
\begin{equation*}
\lim_{\ep \to 0} \E_Z \left[e^{-\theta U^{\ep} (t)} \right]
=
\lim_{\ep \to 0} \exp \left\{ m^{\ep} (t) (e^{-\theta /m^{\ep} (t)} - 1) \right\}
=
e^{- \theta} \,
\end{equation*}
which is our claim for a single $\theta$.

Second, we now extend the above to $\Pbf$-almost surely for {\em all} 
$\theta$.
Let $X^{\ep} (\theta ) = \E_Z \left[e^{-\theta U^{\ep} (t)} \right]$. Then $\{X^{\ep} (\theta) \,;\, \theta \geq 0\}$ is a 
stochastic process whose randomness is based on the random variable $Z$ only.
By \cref{limit.parameter}. 
it is sufficient to prove that $\theta \mapsto X^{\ep} (\theta)$ is a $\Pbf$-almost surely Lipschitz function. To this aim, consider $\theta, h \geq 0$. Then we can express the difference $X(\theta +h) - X(\theta)$ as follows:
\begin{equation*}
X^{\ep} (\theta +h) - X^{\ep} (\theta) = \E_Z \left[ e^{-\theta U^{\ep} (t)} \left( e^{-h U^{\ep} (t)} - 1\right)\right] .
\end{equation*}
We then get
\begin{equation}\label{ineq:diff-X-ep-h-increment}
\left| X^{\ep} (\theta +h) - X^{\ep} (\theta) \right| \leq
\E_Z \left[ \left| e^{-h U^{\ep} (t)} - 1 \right| \right] \leq
h \, \E_Z \left[ U^{\ep} (t) \right] = h \,,
\end{equation}
where the last identity stems from the fact that $\E_Z \left[ U^{\ep} (t) \right] = 1$.
Relation \eqref{ineq:diff-X-ep-h-increment} obviously proves that $\theta \mapsto X^{\ep} (\theta)$ is $1$-Lipschitz $\Pbf$-almost surely and uniformly in $\ep>0$. Then the conclusion follows.

Statement \eqref{NepLimitm0} follows from
\begin{equation}\label{eq:calculate-N-ep/m-ep-beta0}
\frac{N^\ep(t)}{m^\ep(t)} - \frac{\ep^\beta N^\ep(t)}{\bar{m}_0(t)}
= 
\frac{N^\ep(t)}{m^\ep(t)}\left(1-\frac{\ep^\beta m^\ep(t)}{\bar{m}_0(t)}\right)
= 
U^{\ep} (t) \left(
\frac{\bar{m}_0(t) - m^\ep_0(t)}{\bar{m}_0(t)}
\right) .
\end{equation}
Then in the right hand side of \eqref{eq:calculate-N-ep/m-ep-beta0} we have proved that $U^{\ep} (t)$ defined by \eqref{eq:U-ep-fraction} converges to $1$ in $\PP_Z$-probability. Moreover, relation \eqref{m.ep.rate.subcrit} asserts that $\bar m_0 (t) - m^{\ep}_0 (t)$ converges to $0$ $\Pbf$-almost surely. This achieves the proof of \eqref{NepLimitm0}.
\end{proof}

\subsection{Central Limit Theorem: Subcritical Case $\beta < 1$.}
\label{sect:clt-finite-dim-dist}
With the above law of large numbers in hand, we now proceed to investigate
the corresponding central limit theorem. More precisely, in this section,
we will analyze the rescaled process \eqref{centered-Nep} which is copied 
here for convenience:
\begin{equation}
G^\ep(t) =
\ep^{\frac{\beta}{2}}\Big(N^\ep(t) - \frac{\bar{m}_0(t)}{\ep^\beta}\Big).
\tag{\ref{centered-Nep}}
\end{equation}
Here we tackle the subcritical case and establish a
functional limit theorem in the quenched setting.

\begin{theorem}[Functional Central Limit Theorem]\label{thm-FCLT-quenched}
Let $\beta < 2\ka$, where $\ka$ is introduced in \cref{hyp:ergodicity} and $\beta$ is defined in \cref{hyp:arr-times-Gamma}. Then $\Pbf$-almost surely,
\begin{equation}\label{limeq:G-ep-to-bar-G-Pz}
G^\ep(\cdot) \Longrightarrow_{\PP_Z} \bar{G}(\cdot)
\end{equation}
where $\bar{G}$ is a continuous centered Gaussian process. 
The convergence holds in the Skorokhod topology of c\`adl\`ag processes.
\end{theorem}

The proof of \cref{thm-FCLT-quenched} follows by first establishing the 
convergence of marginal and finite dimensional distributions,
in~\cref{prop:CLT-G-ep-1d,prop:CLT-G-ep-md}, respectively. Then the continuity and tightness properties are proved in~\cref{G-continuity} and~\Cref{prop:G-ep-family-tightness}.

\subsubsection{Finite Dimensional Convergence}\label{clt-finite-quenched}
We first establish the finite dimensional convergence statement contained in \cref{thm-FCLT-quenched}. 
We will start with a lemma that will be used in the sequel.

\begin{lemma}\label{lemma:G-ep-EZ-constbound}
Consider $\beta < 2\kappa$ and recall that $\ka \in (0, 1/2)$ is introduced in relation \eqref{eq:hyp-ergodic-on-Z}.
Assume \cref{hyp:ergodicity,hyp:arr-times-Gamma,hyp:F-bar-increments} are fulfilled. Then for any fixed $t \geq 0$ and for any $\ep >0$, we have
\begin{equation}\label{ineq:bound-G-ep-Const}
\E_Z \left[|G^{\ep} (t) | \right] \leq C_t(Z) \,,
\end{equation}
where $C_t (Z)$ is a (random) constant depending on $Z$ only.
\end{lemma}

\begin{proof}
We will make use of properties of Poisson random variable.
To this aim, recall our decomposition \eqref{clt-decompose} for $G^\ep(t)$:
\begin{equation}\label{eq:decomp-G-ep-G1,2-repeat}
G^{\ep} (t)  = G^\ep_1(t) + G^\ep_2(t) \,,
\end{equation}
where $G^{\ep}_1 (t)$ and $G^{\ep}_2 (t)$ are introduced in \eqref{eqdef:Q-ep-1,2-fluctuation}.
By \cref{prop:m-t-of-N-t}, $N^\ep(t)$ is a Poisson random variable with
parameter $m^\ep(t)$. Hence,
\begin{equation*}
\E_Z \left[ |N^{\ep} (t) - m^{\ep} (t) |\right]
\leq
\sqrt{\E_Z \left[\left|N^{\ep} (t) - m^{\ep} (t)\right|^2\right]}
= \sqrt{\text{Var}_Z\left(N^\ep(t)\right)}
= \sqrt{m^\ep(t)} \,.
\end{equation*}
Since one can recast \eqref{eqdef:Q-ep-1,2-fluctuation} as
$G^\ep_1 (t) = \ep^{\beta/2} (N^\ep (t) - m^\ep (t))$, it is easily seen that
\begin{equation}\label{ineq:bound-G-ep1exp}
\E_Z \left[ |G^\ep_1(t) |\right] \leq \sqrt{\ep^\beta m^\ep(t)} 
= \sqrt{m_0^\ep(t)} \,.
\end{equation}
Owing to the fact that the right hand side of \eqref{ineq:bound-G-ep1exp} is uniformly bounded on any interval $[0,T]$, we get that $\E_Z \left[ |G^\ep_1(t) |\right] \leq C_t (Z)$.

For $G^\ep_2(t)$, we have $\E_Z \left[G_2^\ep(t)\right] = \left|G_2^\ep(t) \right|$ as it 
depends on $Z$ only. 
Furthermore, by \eqref{m.ep.rate.subcrit} in \cref{weighted.ergodic}, 
$G^{\ep}_2 (t)$ tends to zero $\Pbf$-almost surely as $\ep\longrightarrow0$.

Combining the above into \eqref{eq:decomp-G-ep-G1,2-repeat}, the conclusion follows.
\end{proof}

\begin{proposition}[Univariate CLT]
\label{prop:CLT-G-ep-1d}
Suppose that \cref{hyp:ergodicity,hyp:arr-times-Gamma,hyp:F-bar-increments} are satisfied and recall that $\ka \in (0, 1/2)$ is introduced in relation \eqref{eq:hyp-ergodic-on-Z}. Then for any fixed $t \geq 0$, we have a quenched central limit theorem for $G^{\ep} (t)$ whenever $0<\beta < 2 \kappa$. More precisely, 
the following holds true $\Pbf$-almost surely,
\begin{equation}\label{eqlim:G-hat-weakconvg}
G^{\ep} (t) \Longrightarrow_{\PP_Z} \mathcal{N} (0,\bar{m}_0(t)) \,.
\end{equation}
\end{proposition}

\begin{proof}
We will make use of the characteristic function of $G^{\ep} (t)$. 
More specifically, we will prove that $\Pbf$-almost surely,
for all $\theta$, the following holds
\begin{equation}\label{eqlim:charfn-G-ep-Gaussian}
\lim_{\ep \to 0} \E_Z \exp\left[i\theta G^\ep(t)\right] = \exp\left[-\frac{\bar{m}_0(t)\theta^2}{2}\right] .
\end{equation}
First, we compute the characteristic function of $G^\ep(t)$ for a fixed $\ep >0$:
\begin{eqnarray}
& & \E_Z \exp\left[i\theta G^\ep(t)\right] \notag\\
& = & \E_Z \exp\left[i\theta \ep^{\frac{\beta}{2}}
\left(N^\ep(t) - \frac{\bar{m}_0(t)}{\ep^\beta}\right)\right]
 =  \exp\left[-i\theta\frac{\bar{m}_0(t)}{\ep^\frac{\beta}{2}}\right]
\E_Z \exp\left[i\theta \ep^{\frac{\beta}{2}} N^\ep(t)\right] .\notag
\end{eqnarray}
Hence, invoking the fact that $N^{\ep} (t)$ is a Poisson random variable with parameter $m^{\ep} (t)$ and owing to considerations similar to \eqref{Laplacetransform:N-ep/m-ep} we get
\begin{eqnarray}
& & \E_Z \exp\left[i\theta G^\ep(t)\right] \notag\\
& = & \exp\left[-i\theta\frac{\bar{m}_0(t)}{\ep^\frac{\beta}{2}}\right]
\exp\left[\frac{m^\ep_0(t)}{\ep^\beta}
\left(e^{i\theta \ep^{\frac{\beta}{2}}}-1\right)\right] \notag \\
& = &
\exp\left[-i\theta\frac{\bar{m}_0(t)}{\ep^\frac{\beta}{2}}\right]
\exp\left[\frac{\bar{m}_0(t)}{\ep^\beta}
\left(e^{i\theta \ep^{\frac{\beta}{2}}}-1\right)\right]
\exp\left[\frac{m^\ep_0(t)-\bar{m}_0(t)}{\ep^\beta}
\left(e^{i\theta \ep^{\frac{\beta}{2}}}-1\right)\right] \notag \\
& = &
\exp\left[\frac{\bar{m}_0(t)}{\ep^\beta}
\left(e^{i\theta \ep^{\frac{\beta}{2}}}-1-i\theta\ep^{\frac{\beta}{2}}\right)\right]
\exp\left[\frac{m^\ep_0(t)-\bar{m}_0(t)}{\ep^\beta}
\left(e^{i\theta \ep^{\frac{\beta}{2}}}-1\right)\right] . \label{eq:char-G-ep-expfactors}
\end{eqnarray}
Now for each fixed $\theta$, we next compute limits for the above expression 
as $\ep\longrightarrow0$.
The convergence of the first exponential factor is done completely in a
deterministic manner,
\begin{equation}\label{lim:firstexpofactor}
\lim_{\ep\to0}\exp\left[\frac{\bar{m}_0(t)}{\ep^\beta}
\left(e^{i\theta \ep^{\frac{\beta}{2}}}-1-i\theta\ep^{\frac{\beta}{2}}\right)\right]
= \exp\left[-\frac{\bar{m}_0(t)\theta^2}{2}\right].
\end{equation}
For the convergence of the second exponential factor in 
\eqref{eq:char-G-ep-expfactors}, we rewrite it as
\begin{equation*}
\exp\left[\frac{m^\ep_0(t)-\bar{m}_0(t)}{\ep^\beta}
\left(e^{i\theta \ep^{\frac{\beta}{2}}}-1\right)\right]
=\exp\left[\left(\frac{m^\ep_0(t)-\bar{m}_0(t)}{\ep^{\beta/2}} \right)
\left(\frac{e^{i\theta \ep^{\frac{\beta}{2}}}-1}{\ep^{\beta/2}}\right) \right] .
\end{equation*}
Then by \eqref{m.ep.rate.subcrit} of \cref{weighted.ergodic} and the fact that
$\lim_{\ep \to 0} |\frac{e^{i\theta \ep^{\beta /2}}-1}{\ep^{\beta/2}}| = \theta$,
the convergence of the above expression follows: $\Pbf$-almost surely,
\begin{equation}\label{lim:secexpofactor}
\lim_{\ep\to0}\exp\left[\frac{m^\ep_0(t)-\bar{m}_0(t)}{\ep^\beta}
\left(e^{i\theta \ep^{\frac{\beta}{2}}}-1\right)\right]
=1 \,.
\end{equation}
Hence, plugging \eqref{lim:firstexpofactor} and \eqref{lim:secexpofactor} into \eqref{eq:char-G-ep-expfactors}, we get that \eqref{eqlim:charfn-G-ep-Gaussian} is valid for each fixed $\theta$.

Next we will establish the convergence $\Pbf$-almost surely for all $\theta$.
This follows from \cref{limit.parameter} upon showing that the following
map is Lipschitz:
\begin{equation}\label{map:theta-to-X-ep-Lipschitz}
\theta \longmapsto X^\ep(\theta) 
:= \mathbb{E}_Z\left[\exp\left(i\theta G^\ep(t)\right)\right].
\end{equation}
To prove this, we compute
\begin{eqnarray*}
\left|X^{\ep} (\theta +h) - X^{\ep} (\theta)\right| 
& = & 
\left|\E_Z \left[ e^{i\theta G^{\ep} (t)} \left( e^{i h G^{\ep} (t)} - 1\right)\right]\right|
\\
& \leq & \E_Z \left|e^{i h G^{\ep} (t)} - 1\right|\\
& \leq & \E_Z \left|hG^{\ep} (t)\right|\\
& \leq & C_t (Z) \,h \,,
\end{eqnarray*}
where we have used \cref{lemma:G-ep-EZ-constbound} for the last inequality.
The above upper bound obviously shows that $\theta \longmapsto X^\ep(\theta)$ in \eqref{map:theta-to-X-ep-Lipschitz} is Lipschitz, which
concludes the proof of the whole proposition.
\end{proof}

The second step toward the proof of \cref{thm-FCLT-quenched} is the 
following multivariate CLT.
\begin{proposition}[Multivariate CLT]
\label{prop:CLT-G-ep-md}
Let the vector $(G^{\ep} (t_1), G^{\ep} (t_2), \ldots , G^{\ep} (t_n))^T$ be defined according to \eqref{centered-Nep}.
Then $\Pbf$-almost surely, for any $0 \leq t_1 < t_2 < \cdots < t_n$ and
as $\ep \to 0$, the vector
weakly converges to a centered Gaussian vector:
\begin{equation}\label{lim:mult-G-ep-t1-tn}
\left(G^{\ep} (t_1) , G^{\ep} (t_2) , \ldots , G^{\ep} (t_n) \right)^T
\xrightarrow{(\dd)}
\mathcal{N} (0, \Gamma) \,,
\end{equation}
where $\Gamma \in \R^{n \times n}$ denotes the covariance matrix such that $\Gamma = \left\{\Gamma (t_i , t_j) \,;\, i,j \in \{1,2, \ldots , n\}\right\}$ with
\begin{equation}\label{defeq:Gamma-titj-inmultdim}
\Gamma (t_i , t_j) = \bar\psi \int_{0}^{t_i \wedge t_j} \la (s) \bar{F}_s (t_i \vee t_j - s) \,\dd s \,.
\end{equation}
\end{proposition}




The proof of the above will be presented for the interested reader in \cref{sect:appendix-prop-multclt}.
We will concentrate on the bivariate case as it already 
illustrates
the main idea and computation underlying the proof. 
As a preparation for this, we will
introduce some notations.
In the sequel, let $0< t_1< t_2$ be two fixed times.
\begin{enumerate}
\item
Inside the quadrant $\R_+ \times \R_+$, we consider the
three disjoint regions $\{A_i \,, \, i = 1,2,3\}$ defined as follows (see \cref{fig:A_i-region-add} for a similar construction):
\begin{eqnarray}
A_1 &=& \{(\ga, l) \in \mathbb{R}_+ \times \mathbb{R}_+ : \ga \leq t_1  \text{ and }  t_1 < \ga+l \leq t_2\} \,;
\label{def:region_A1}\\
A_2 &=& \{(\ga, l)\in \mathbb{R}_+ \times \mathbb{R}_+ : \ga \leq t_1   \text{ and }  t_2 < \ga + l\} \,;
\label{def:region_A2}\\
A_3 &=& \{(\ga, l) \in \mathbb{R}_+ \times \mathbb{R}_+ : t_1 < \ga \leq t_2   \text{ and }  t_2 < \ga+l\} \,.
\label{def:region_A3}
\end{eqnarray}

\item
Since the $A_i$’s are disjoint, the quantities $\{M^{\ep} (A_i); \, i = 1,2,3\}$ are independent Poisson random variables under the quenched probability $\PP_Z$. Following the proofs of \cref{prop:m-t-of-N-t,prop:gen-limit-m-ep}, we deduce that 
their respective quenched means are given by
\begin{equation}\label{eq:def-mean-M-ep-A_i-reg}
m^{\ep}(A_i) := \mathbb{E}_Z [M^{\ep} (A_i)] = \frac{1}{\ep^{\beta}} \int_{A_i} \nu (s, \dd r) \la (s) \psi (Z_{s/\ep}) \dd s \,, \quad \text{for}\,\,\,i=1,2,3 \,,
\end{equation}
with limits
\begin{eqnarray}
\bar{m}_0(A_1) = \lim_{\ep \to 0} \ep^{\beta} m^{\ep}(A_1)
&=& \lim_{\ep \to 0} \ep^{\beta} \mathbb{E}_Z [M^{\ep} (A_1)] = (\Lambda (t_1) - \Lambda (t_2)+ \Xi_{t_1t_2}) \bar{\psi}
\notag \\
\bar{m}_0(A_2) = \lim_{\ep \to 0} \ep^{\beta} m^{\ep}(A_2)
&=& \lim_{\ep \to 0} \ep^{\beta} \mathbb{E}_Z [M^{\ep} (A_2)] = (\Lambda (t_2) - \Xi_{t_1t_2}) \bar{\psi}
\label{eqns:m-ep-i-limit} \\
\bar{m}_0(A_3) = \lim_{\ep \to 0} \ep^{\beta} m^{\ep}(A_3)
&=& \lim_{\ep \to 0} \ep^{\beta} \mathbb{E}_Z [M^{\ep} (A_3)] = \Xi_{t_1t_2} \bar{\psi} ,
\notag
\end{eqnarray}
where $\Lambda (t)$ and $\bar\psi$ are defined previously by 
\eqref{eq:quantity-sigma-psi-bar} and $\Xi_{t_1t_2}$ is equal to
\begin{equation}\label{eq:def-Xi-t1t2}
\Xi_{t_1t_2} = \int_{t_1}^{t_2} \la (s) \bar{F}_s (t_2 -s) \, \dd s \,.
\end{equation}

\item
Next we consider the following rescaled version of $M^{\ep}(A_i)$,
\begin{equation}\label{defeq:M-ep-A-i-m-ep-i}
\widetilde{M^\ep}(A_i) 
= \ep^{\beta/2} \left(M^{\ep} (A_i) - m^{\ep}(A_i) \right),
\quad \text{for} \quad i = 1,2,3 .
\end{equation}
Based on the univariate CLT 
(\cref{prop:CLT-G-ep-1d}) and its proof,
the following quenched CLT holds similarly: $\Pbf$-almost surely,
\begin{equation}\label{eq:note-for-recaled-M-ep}
\widetilde{M^\ep}(A_i) 
\xrightarrow{(\dd)} \mathcal{N} (0,\bar m_0(A_i)), 
\quad \text{for} \quad i = 1,2,3.
\end{equation}
\end{enumerate}

\noindent
With the above, 
we now state our main central limit result for bivariate vectors.

\begin{proposition}[Bivariate CLT]\label{prop:bivariate-clt-for-G-ep}
Let the vector $(G^{\ep} (t_1), G^{\ep} (t_2) )$ be defined according to the expression in ~\eqref{centered-Nep}. 
Then as $\ep \to 0$ and $\Pbf$-almost surely, $(G^{\ep} (t_1), G^{\ep} (t_2) )$ weakly converges to a centered Gaussian process:
\begin{equation}\label{eq:clt-for-G-ep-vec}
\begin{pmatrix}
 G^{\ep} (t_1) \\
 G^{\ep} (t_2)
\end{pmatrix}
\xrightarrow{(\dd)} \mathcal{N} (0, \Gamma) ,
\quad \text{with} \quad
\Gamma = \begin{pmatrix} 
\bar m_0(A_1) + \bar m_0(A_2) & \bar m_0(A_2)\\ \bar m_0(A_2) & \bar m_0(A_2)+\bar m_0(A_3) \end{pmatrix} .
\end{equation}
As in \cref{prop:CLT-G-ep-1d} relation \eqref{eqlim:G-hat-weakconvg}, 
the convergence in \eqref{eq:clt-for-G-ep-vec} holds in distribution with respect to the quenched probability $\PP_Z$.
\end{proposition}

\begin{proof}
We observe that the vector $(G^{\ep} (t_1) , G^{\ep} (t_2))$ can be decomposed as
\begin{equation}\label{eq:G-ep-decompmatrix-M-ep}
\begin{pmatrix}
 G^{\ep} (t_1) \\
 G^{\ep} (t_2)
\end{pmatrix}
= \begin{pmatrix}1 & 1 & 0 \\ 0 & 1 & 1\end{pmatrix}
\begin{pmatrix} 
\widetilde{M^{\ep}}(A_1) \\ 
\widetilde{M^{\ep}}(A_2) \\ 
\widetilde{M^{\ep}}(A_3) 
\end{pmatrix} ,
\end{equation}
where the random variables $\widetilde{M^{\ep}}(A_i)$’s are introduced in 
\eqref{defeq:M-ep-A-i-m-ep-i}. 
Hence with the help of \eqref{eq:note-for-recaled-M-ep}, together with the fact that $\widetilde{M^{\ep}}(A_i)$’s are independent, 
we get the following multivariate quenched limit in distribution as $\ep \to 0$:
\begin{equation*}
\begin{pmatrix} 
\widetilde{M^{\ep}}(A_1) \\ 
\widetilde{M^{\ep}}(A_2) \\ 
\widetilde{M^{\ep}}(A_3) 
\end{pmatrix}
\xrightarrow{(\dd)}
\mathcal{N} \left(0, I_{\bar m} \right) \quad \text{with} \quad I_{\bar m} =
\begin{pmatrix} \bar m_0(A_1) & 0 & 0\\ 0& \bar m_0(A_2) &0\\ 0&0& \bar m_0(A_3)  \end{pmatrix} .
\end{equation*}
Therefore, some standard considerations for Gaussian vectors yield:
\begin{equation*}
\begin{pmatrix}
 G^{\ep} (t_1) \\
 G^{\ep} (t_2)
\end{pmatrix}
= \begin{pmatrix}1 & 1 & 0 \\ 0 & 1 & 1\end{pmatrix}
\begin{pmatrix} 
\widetilde{M^{\ep}}(A_1) \\ 
\widetilde{M^{\ep}}(A_2) \\ 
\widetilde{M^{\ep}}(A_3) 
\end{pmatrix}
\xrightarrow{(\dd)} \mathcal{N} (0, \Gamma) \,,
\end{equation*}
where the covariance matrix $\Gamma$ is computed as follows:
\begin{equation*}
\Gamma = \begin{pmatrix}1 & 1 & 0 \\ 0 & 1 & 1\end{pmatrix} I_{\bar m}
\begin{pmatrix} 1 & 0\\ 1&1 \\ 0& 1 \end{pmatrix}
= \begin{pmatrix} \bar m_0(A_1) + \bar m_0(A_2) & \bar m_0(A_2)\\ \bar m_0(A_2) & \bar m_0(A_2)+\bar m_0(A_3) \end{pmatrix} .
\end{equation*}
Therefore, we may conclude.
\end{proof}

\begin{remark}\label{rmk:expression-Gamma-titj}
We point out here that formula \eqref{eq:clt-for-G-ep-vec} is actually 
consistent with \eqref{defeq:Gamma-titj-inmultdim}.
Indeed, combining \eqref{eqns:m-ep-i-limit}, \eqref{eq:def-Xi-t1t2}, and 
\eqref{eq:quantity-sigma-psi-bar}, one can readily check that 
\begin{equation*}
\E \left[ G(t_1) G(t_2) \right] = \bar m_0(A_2) =
\bar\psi \int_{0}^{t_1} \la (s) \bar{F}_s (t_2 -s) \,\dd s \,.
\end{equation*}
Therefore, the matrix $\Gamma$ in \eqref{eq:clt-for-G-ep-vec} can be written as $\Gamma = \left\{\Gamma (t_i , t_j) \,;\, i,j \in \{1,2\}\right\}$ with
\begin{equation}\label{eq:gamma-titj-with-psi-bar}
\Gamma (t_i , t_j) = \bar\psi \int_{0}^{t_i \wedge t_j} \la (s) \bar{F}_s (t_i \vee t_j - s) \,\dd s \,.
\end{equation}
\end{remark}

The above univariate and bivariate results already reveal 
some interesting features. The following 
proposition partially establishes the ``fractional Brownian motion type
self-similarity’’ property for the limit of $G^\ep(t)$.

\begin{proposition}\label{prop:self-similar-process-fbm}
For $t>0$, let $G(t)$ be a random variable distributed according to 
\eqref{eqlim:G-hat-weakconvg}, that is, $G^{\ep} (t) \sim \mathcal{N} (0,\bar{m}_0(t))$.
We also assume some non-degeneracy conditions on $\la$ and $\psi$, namely, the existence of $\underline\la$ and $\overline\la$ such that
\begin{equation}\label{eq:nondegeneracy-condition-la-psi-1}
0<\underline\la \leq \la (s) \leq \overline\la < \infty \,, \quad \text{and} \quad 0< \bar\psi < \infty \,.
\end{equation}
Then there exist two constants $0 < c_1 <c_2< \infty$ such that for all $t \geq 1$ we have
\begin{equation}\label{ineq:expected-G(t)}
\begin{dcases}
c_1 \, t^{1-\al} \,\, \leq \,\, \E \left[G(t)^2\right] \,\,\leq \,\, c_2 \, t^{1-\al} \,, & \text{if} \quad \al \in (0,1) \\
c_1 \,\, \leq \,\, \E \left[G(t)^2\right] \,\,\leq \,\, c_2 \,, & \text{if} \quad \al \in (1, \infty)
\end{dcases}
\,.
\end{equation}
In particular, for $\al \in (0,1)$, $G$ has an fBm-type behavior, with Hurst parameter $H = \frac{1}{2} (1-\al) \in (0, \frac{1}{2})$ while for $\al \in (1, \infty)$, $G$ has a stationary-noise-type behavior.
\end{proposition}

\begin{proof}
Recall from \eqref{eq:limit-epbeta-m-ep} that we have $\bar{m}_0 (t) = \Lambda (t) \bar{\psi}$ and we have assumed $\bar\psi >0$. Therefore, our claim \eqref{ineq:expected-G(t)} is reduced to proving the following for $t \geq 1$,
\begin{equation}\label{ineq:sigma-t-bounds-claim}
\begin{cases}
c_1 \, t^{1-\al} \,\, \leq \,\, \Lambda (t) \,\, \leq \,\, c_2 \, t^{1-\al} \,, & \text{if} \quad \al \in (0,1) \\
c_1 \,\, \leq \,\, \Lambda (t) \,\,\leq \,\, c_2 \,, & \text{if} \quad \al \in (1, \infty)
\end{cases}
\,,
\end{equation}
where $\Lambda (t)$ is defined by \eqref{eq:quantity-sigma-psi-bar}. Moreover, thanks to \eqref{def:F-bar-generalization} and \eqref{eq:nondegeneracy-condition-la-psi-1}, one can write
\begin{equation}\label{ineq:sigma-t-bounds-cK}
c_1 \, \underline\la \iot \left( \frac{1}{(t-s)^\alpha} \wedge 1 \right) \dd s \,\, \leq \,\,
\Lambda (t) \,\, \leq \,\,
c_2 \, \overline\la \iot \left( \frac{1}{(t-s)^\alpha} \wedge 1 \right) \dd s \,.
\end{equation}
Starting from \eqref{ineq:sigma-t-bounds-cK}, explicit calculations
yield the desired result in \eqref{ineq:sigma-t-bounds-claim}, completing 
our proof.
\end{proof}

We now state a Brownian type regularity result for the process $\bar G$.

\begin{proposition}\label{G-continuity}
Let $\bar{G}(\cdot)$ be the centered Gaussian process given in the limit \eqref{limeq:G-ep-to-bar-G-Pz} obtained in \cref{thm-FCLT-quenched}.
Then $\bar{G}(\cdot)$ is a continuous stochastic process.
More specifically, $\bar G$ is a $\mu$-H\"older continuous process on any interval $[0,T]$, for all $\mu \in (0, 1/2)$.
\end{proposition}

\begin{proof}
Since $\bar G$ is a Gaussian process, our claim is reduced to showing that for $s < t$,
\begin{equation}\label{defeq:R-st-E-Gst}
R_{st} \equiv \E \left[(\bar G (t) - \bar G (s))^2\right] \,\lesssim\, t-s \,.
\end{equation}
(As $\bar G(\cdot)$ is Gaussian, the above leads to
that $\E \left[|\bar G (t) - \bar G (s)|^{2m}\right] \,\lesssim\, |t-s|^m$
for all $m > 1$. The continuity of $\bar G(\cdot)$ then follows
from the Kolmogorov continuity theorem.)
Recall that the covariance function $\Gamma$ of $\bar G$ is given by \eqref{defeq:Gamma-titj-inmultdim}. Furthermore, an elementary algebraic manipulation reveals that for $s< t$, we have
\begin{equation*}
R_{st} = \Gamma (t,t) - 2 \Gamma (t,s) + \Gamma (s,s) \,.
\end{equation*}
This rectangular increment can be recast as
\begin{equation*}
R_{st} = \bar\psi \int_{s}^{t} \la (u) \bar F_u (t-u) \,\dd u
+ \bar\psi \ios \la (u) \left(\bar F_u (s-u) - \bar F_u (t-u)\right) \dd u \,.
\end{equation*}
Furthermore, invoking \eqref{def:F-bar-generalization}, \eqref{ineq:partial-ell-s-bounds} and the fact that $\la$ is a bounded function (see \cref{hyp:arr-times-Gamma}), we get that
\begin{equation*}
\left|R_{st}\right|  \leq c_{\bar\psi, \,\la} \,(t-s) \,.
\end{equation*}
Our claim \eqref{defeq:R-st-E-Gst} is thus verified, which completes the proof.
\end{proof}

\subsubsection{Tightness}\label{clt-tightness-quenched}
The final step to establish
\cref{thm-FCLT-quenched}
is proving the tightness of the sequence of processes $\{G^{\ep} \,;\, \ep >0\}$.
We will establish this in the space of 
c\`adl\`ag functions, $D\left([0, \infty)\right)$. 
However, as stated in \cref{G-continuity}, the limit process $\bar{G}$ 
in fact lives in the space of continuous functions.

\begin{proposition}\label{prop:G-ep-family-tightness}
Let $\{G^{\ep} \,;\, \ep >0\}$ be the family of processes such that $G^{\ep} (t)$ is defined by \eqref{centered-Nep}. Then $\Pbf$-almost surely, the sequence $\{G^{\ep} \,;\, \ep >0\}$ is tight in the space of c\`adl\`ag functions $D$ equipped with the Skorokhod metric.
\end{proposition}

\begin{proof}
We use a chaining type argument which is inspired by \cite{RR} but has to be adapted to our quenched setting. First, notice that the tightness of the sequence of processes $\{G^{\ep} \,;\, \ep >0\}$ for the Skorohod topology in $D\left([0, \infty)\right)$ is easily deduced from the tightness in $D\left([0, K_0]\right)$ for an arbitrarily large $K_0 > 0$. For simplicity of notation and without loss of generality, we will let $K_0 = 1$ in the sequel.

First, it is shown in \cite[Theorem 15.5]{Bill} that the tightness is obtained by controlling the increments of $G^{\ep}$, $\Pbf$-almost surely. Namely, it suffices to prove that for all rational $\theta , \eta >0$, there exists a set $\Omega_{\theta, \eta}$ such that $\Pbf (\Omega_{\theta, \eta}) = 1$ and for all $\omega \in \Omega_{\theta, \eta}$, there exist $\ep_0 (\omega) , \delta (\omega) >0$ such that for any $\ep \leq \ep_0 (\omega)$ we have
\begin{equation}\label{eqn3.12}
\PP_{Z(\omega)} \left[ \sup_{|t_2-t_1| \leq \delta (\omega) \,;\, 0 \leq t_1 < t_2 \leq 1} |G^{\ep} (t_2) - G^{\ep} (t_1)| > \theta\right]
< \eta \,. 
\end{equation}

Next, we discretize the parameter $\delta$ involved in the supremum norm in \eqref{eqn3.12}. Specifically, one can choose $i = i (\omega)$ and $v_0 = v_0 (\ep_0 (\omega))$ 
such that for $\delta = 2^{-i}$ and $v \geq v_0$, we have
\begin{multline}\label{ineq:sup-increments-G-ep-t1t2}
\sup_{|t_2-t_1| < \delta} |G^{\ep} (t_2) - G^{\ep} (t_1)| \leq \\
\max_{t_1, t_2 \in F_v , |t_2-t_1| \leq \delta} \left|G^{\ep} (t_2) - G^{\ep} (t_1)\right| +
2 \bigvee_{k =0}^{2^v-1} \sup_{0 \leq u \leq 2^{-v}} \left|G^{\ep} (k2^{-v} +u) - G^{\ep} (k 2^{-v})\right| ,
\end{multline}
where the discrete set $F_v$ is defined by
\begin{equation}\label{def:F_j-dyadic}
F_v := \left\{k 2^{-v} : 0 \leq k \leq 2^v \right\} .
\end{equation}

Putting together \eqref{eqn3.12} and \eqref{ineq:sup-increments-G-ep-t1t2}, the desired tightness property is thus reduced to the following two statements:

\begin{enumerate}
\item For all rational $\theta ,  \eta >0$, one can find a set $\Omega_{\theta, \eta}$ with $\Pbf (\Omega_{\theta, \eta}) = 1$ and for all $\omega \in \Omega_{\theta, \eta}$, there exist $\ep_0 (\omega),  \delta (\omega) >0$ satisfying the following: for all $\ep \leq \ep_0 (\omega)$, there exists $v_0 = v_0 (\ep_0 (\omega))$ such that for $v \geq v_0$ we have
\begin{equation}\label{ineq:prob-sup-diff-G-ep}
\PP_{Z (\omega)} \left[ \sup_{t_1, t_2 \in F_v \,;\, |t_2-t_1| \leq \delta (\omega)} |G^{\ep} (t_2) - G^{\ep} (t_1)| > \theta\right]
< \eta \,,
\end{equation}
where $F_v$ is defined by \eqref{def:F_j-dyadic}.

\item As $\ep \to 0$ and $v = v (\ep) \to \infty$, the following limit holds true in $\PP_Z$-probability,
\begin{equation}\label{lim:max-increment-G-ep-v}
\bigvee_{k =0}^{2^v-1} \sup_{0 \leq u \leq 2^{-v}} |G^{\ep} (k2^{-v} +u) - G^{\ep} (k 2^{-v})| \to 0 \,.
\end{equation}
\end{enumerate}

Relations \eqref{ineq:prob-sup-diff-G-ep} and \eqref{lim:max-increment-G-ep-v} are respectively proved in \cref{lemma:tightness-ineqstate1,lemma:tightness-ineqstate2}. Thus, we conclude the proof of the tightness property.
\end{proof}

\section{Annealed Analysis}\label{Sect:annealed-analysis}
In this section we investigate $N^\ep(\cdot)$ in the annealed regime, 
i.e. with respect to the probability measure $\PP$ and $\E$. 
To illustrate the idea behind our analysis, we recall the notation given in \eqref{def:anneal-P}. 
Specifically, for any event $A\subset\Omega$ and random variable $X$, we have
$$
\PP(A) = \Pbf\big(\PP_Z(A)\big),
\quad\text{and}\quad
\E [X] = \mathbf{E}\big[\E_ZX\big].
$$
The previous two sections concerned the expressions $\PP_Z(A)$ and $\E_Z [X]$,
which depend on each realization of $Z$. Hence to extend the
result to $\PP(A)$ and $\E [X]$, we need some uniformity or regularity 
with respect to $Z$. This is captured by \cref{hyp:natual-annealed-bound} 
(which is an $L^1$-ergodicity condition)
in replacement of \cref{hyp:ergodicity}
(which is an almost sure condition). In particular, the $L^1$-ergodicity
implies the existence of a subsequence as $t_i\longrightarrow\infty$ so that
for $\Pbf$-almost every realization of $Z$, it holds that
\begin{equation*}
\frac{1}{t_i}\int_0^{t_i}\psi(Z_s)\,\dd s\longrightarrow\bar{\psi}.
\end{equation*}
Therefore, we will partially be able to take advantage of the quenched results obtained in \cref{Sect:LLN-for-N-ep,sect:clt-finite-dim-dist}. This will be
made precise in the actual proof of the statements below.

\subsection{Law of Large Numbers for $N^\ep(t)$}
As in the quenched case, the following statement is 
applicable for both the subcritical and supercritical regimes.
\begin{theorem}[Annealed LLN]\label{thm-LLN-annealed}
Let $N^{\ep}$ be the process defined by \eqref{eq:N-t1}, and assume that 
\cref{hyp:arr-times-Gamma,hyp:F-bar-increments,hyp:natual-annealed-bound} 
are satisfied. Furthermore, suppose $\beta > 0$. Then for any fixed $t \geq 0$, we have
\begin{equation}\label{eq:LLN-N-ep-equal-1-ann}
\lim_{\ep \to 0} \frac{N^{\ep} (t)}{m^{\ep} (t)} = 1
\quad \text{in $\PP$-probability},
\end{equation}
where $m^{\ep} (t)$ is defined by \eqref{eq:m-t-def}.
More precisely, one can restate the above as follows:
$\PP$-almost surely we have that for all $\delta > 0$,
\begin{equation}\label{limeq:P-lln-m-ep-delann}
\lim_{\ep \to 0} \, \PP\left(\left|\frac{N^{\ep} (t)}{m^{\ep} (t)} - 1 \right| \geq \delta \right) = 0 \,.
\end{equation}
Furthermore, we have
\begin{equation}\label{limeq:lln-epN-ann}
\lim_{\ep \to 0} \frac{\ep^{\beta} N^{\ep} (t)}{\bar{m}_0(t)} 
= 1 \,, \quad \text{in $\PP$-probability} \,.
\end{equation}
\end{theorem}

\begin{proof}
Similar to its counterpart in \cref{thm-LLN-quenched} for the quenched 
analysis, we define
\begin{equation}\label{eq:U-ep-fractionann}
U^{\ep} (t) = \frac{N^{\ep} (t)}{m^{\ep} (t)} \,,
\end{equation}
and we will prove that $\lim_{\ep \to 0} U^{\ep} (t) = 1$ in $\PP$-probability. Since $1$ is a constant, we will again make use of the well-known fact that 
convergence in probability stems from a mere convergence in distribution. 
For this, we compute Laplace transform of $U^\ep(t)$ and 
prove that for all $\theta \geq 0$:
\begin{equation}\label{eq:limit-exp-theta-U-epann}
\lim_{\ep \to 0}
\E\left[e^{-\theta U^{\ep} (t)} \right]
= \lim_{\ep \to 0} \mathbf{E} \Big\{\E_Z\left[e^{-\theta U^{\ep} (t)} \right]\Big\}
= e^{-\theta}.
\end{equation}

Note that by \cref{hyp:natual-annealed-bound} and
following the proof of \cref{thm-LLN-quenched}, 
we have {\em up to a subsequence $\ep_i\rightarrow0$}, for 
$\Pbf$-almost surely that
$\E_Z\left[e^{-\theta U^{\ep_i} (t)} \right]$ converges to $e^{-\theta}$.
As $\E_Z\left[e^{-\theta U^{\ep_i} (t)} \right]$ is bounded by one, by
the Lebesgue Dominated Convergence Theorem, we also have for the same
subsequence that
$\mathbf{E} \Big\{\E_Z\left[e^{-\theta U^{\ep_i} (t)} \right]\Big\}$
converges to $e^{-\theta}$. As the limit is uniquely defined, 
i.e. it does not
depend on the subsequence, the convergence can be extended to the
whole sequence $\ep\rightarrow0$. This achieves the proof of \eqref{limeq:P-lln-m-ep-delann}.

The proof of \eqref{limeq:lln-epN-ann} is left to the patient reader. It is based on the corresponding statement in \cref{thm-LLN-quenched}, plus Lebesgue Dominated Convergence Theorem (which can be applied since $\psi$ is bounded). 

We can thus conclude out proof.
\end{proof}

\subsection{Central Limit Theorem: Subcritical Case $\beta < 1$}
For this subsection, the notation $G^\ep(\cdot)$ still refers to the
rescaled process defined by \eqref{centered-Nep}.
As mentioned in \cref{hyp:ergodicity}, we have assumed that $\ka$ in \eqref{eq:hyp-ergodic-on-Z} could be taken as close as we wish to $1/2$. Hence, for this section we will assume (without loss of generality) that $\beta < 2 \ka <1$. Our main annealed CLT type result in this context can be summarized as follows.

\begin{theorem}[Functional CLT, Subcritical, Annealed]
\label{thm-FCLT-sub-annealed}
\begin{equation}
G^{\ep} (\cdot) \Longrightarrow_{\PP} \bar{G}(\cdot)
\end{equation}
where $\bar{G}(\cdot)$ is the continuous centered Gaussian process featured in \cref{thm-FCLT-quenched}, with variance and covariance functions respectively introduced in \eqref{eqlim:G-hat-weakconvg} and \eqref{defeq:Gamma-titj-inmultdim}.
The convergence holds in the Skorokhod topology of c\`{a}dl\`{a}g processes.
\end{theorem}

The proof involves first establishing convergence of the marginal and finite dimensional distributions in~\Cref{prop:CLT-G-ep1d-annealed,prop:multivariate-clt-annealed} and then proving tightness in~\Cref{prop:tightness-annclt}. The theorem is thus proved by appealing to~\cite[Theorem 15.1]{Bill}. Note that this result does {\it not} follow from~\cite{RR} since the results there crucially depend on the presumption that the process $M^\epsilon$ is a Poisson random measure; however, this is clearly not the case under the annealed measure, since the performance measures must be averaged over the random environment process. 

We start by establishing the convergence of the marginal distribution under the annealed measure.

\begin{proposition}[Univariate CLT]\label{prop:CLT-G-ep1d-annealed}
Suppose that 
\cref{hyp:arr-times-Gamma,hyp:F-bar-increments,hyp:natual-annealed-bound} 
are satisfied.
Then for any fixed $t \geq 0$, we have an annealed central limit theorem for $G^{\ep} (t)$ whenever $0<\beta < 2 \kappa$. More precisely,
\begin{equation}\label{eqlim:G-hat-weakconvg-ann}
G^{\ep} (t) \Longrightarrow_{\PP} \mathcal{N} (0,\bar{m}_0(t)) \,,
\end{equation}
where we recall from \cref{prop:gen-limit-m-ep} that $\bar{m}_0(t)$ is defined by \eqref{eq:lim-m-ep}.
\end{proposition}

\begin{proof}
As in \cref{prop:CLT-G-ep-1d} of quenched analysis, we will prove this by making use of characteristic function of $G^\ep(t)$. More specifically, it suffices to prove that
\begin{equation}\label{eqlim:charfn-G-ep-Gaussianann}
\lim_{\ep \to 0} \E \exp\left[i\theta G^\ep(t)\right] =
\lim_{\ep \to 0} \mathbf{E} \Big\{\E_Z \exp\left[i\theta G^\ep(t)\right]\Big\}
= \exp\left[-\frac{\bar{m}_0(t)\theta^2}{2}\right] ,
\end{equation}
from which our claim \eqref{eqlim:G-hat-weakconvg-ann} immediately follows.

Again, by \cref{hyp:natual-annealed-bound}, following the proof of
\cref{prop:CLT-G-ep-1d}, we can show that 
$\Pbf$-almost surely, for a subsequence $\ep_i\rightarrow0$, 
the following holds true:
\begin{equation*}
\lim_{\ep_i \to 0} \E_Z \exp\left[i\theta G^{\ep_i}(t)\right] 
= \exp\left[-\frac{\bar{m}_0(t)\theta^2}{2}\right] .
\end{equation*}
Furthermore, as $\left|\exp [i\theta G^\ep(t)] \right|\leq 1$, we have
$\Big| \E_Z \exp\left[i\theta G^\ep(t)\right] \Big|\leq 1$.
Thus by Lebesgue’s dominated convergence theorem, our claim 
\eqref{eqlim:charfn-G-ep-Gaussianann} is achieved for the subsequence
$\ep_i\rightarrow0$. As the limit is again uniquely defined, 
the convergence can be extended to the whole sequence $\ep\rightarrow0$.
This completes the proof of the proposition.
\end{proof}

Next we will state the multivariate CLT for 
multivariate random variables of the form $(G^{\ep} (t_1), G^{\ep} (t_2), \ldots , G^{\ep} (t_n))$ with fixed $t_0=0 < t_1< t_2< \cdots <t_n$.
\begin{proposition}[Multivariate CLT]
\label{prop:multivariate-clt-annealed}
Let \cref{hyp:arr-times-Gamma,hyp:F-bar-increments,hyp:natual-annealed-bound} prevail.
As $\ep \to 0$, we have that
$(G^{\ep} (t_1), G^{\ep} (t_2), \ldots , G^{\ep} (t_n))^T$ weakly converges to a centered Gaussian vector:
\begin{equation}\label{lim:mult-G-ep-t1-tn-ann}
\left(G^{\ep} (t_1) , G^{\ep} (t_2) , \ldots , G^{\ep} (t_n) \right)^T
\Longrightarrow_{\PP}
\mathcal{N} (0, \Gamma)
\end{equation}
where $\Gamma \in \R^{n \times n}$ denotes the covariance matrix such that $\Gamma = \left\{\Gamma (t_i , t_j) \,;\, i,j \in \{1,2, \ldots , n\}\right\}$ with
\begin{equation}\label{defeq:Gamma-titj-inmultdim-ann}
\Gamma (t_i , t_j) = \bar\psi \int_{0}^{t_i \wedge t_j} \la (s) \bar{F}_s (t_i \vee t_j - s) \,\dd s \,.
\end{equation}
\end{proposition}

\begin{proof}
This proof combines the proof of \cref{prop:CLT-G-ep-md} (based on the decomposition using the sets $A_i$) and the dominated convergence arguments from \cref{prop:CLT-G-ep1d-annealed}. Details are left to the interested reader.
\end{proof}


To establish our functional central limit theorem in the annealed regime, we shall prove a tightness property for the family of processes $\{G^{\ep} \,;\, \ep >0\}$.
This property will rely on some further assumptions concerning the coefficients of our queueing system.

\begin{hypothesis}\label{supercri-hyp:coef-la-psi}
Let $\la$ and $\psi$ be the functions featuring in \cref{hyp:arr-times-Gamma}. We assume that there exist four constants $\underline\la \,, \overline\la >0$ and $\underline\psi \,, \overline\psi >0$ such that
\begin{equation}\label{ineqs:la-psi-ulbdd}
\underline\la \leq \la (t) \leq \overline\la \,, \qquad \underline\psi \leq \psi (z) \leq \overline\psi \,, \quad
\text{for all} \quad t \in \R_+ \,, z \in \R^d \,.
\end{equation}
Notice that the hypothesis on $\la$ in \eqref{ineqs:la-psi-ulbdd} has already been used in \cref{prop:self-similar-process-fbm}. 
\end{hypothesis}

With this new hypothesis in hand, the tightness result can be spelled out as follows.

\begin{proposition}[{Tightness for CLT in the Annealed Regime}]\label{prop:tightness-annclt}
Assume \cref{hyp:arr-times-Gamma,hyp:F-bar-increments,hyp:natual-annealed-bound,supercri-hyp:coef-la-psi} hold true.
Let $\{G^{\ep} \,;\, \ep >0\}$ be the family of processes such that $G^{\ep} (t)$ is defined by \eqref{centered-Nep}. Then in $\PP$-probability, the sequence $\{G^{\ep} \,;\, \ep >0\}$ is tight in the space of c\`adl\`ag functions $D ([0, \infty))$ equipped with the Skorokhod metric.
\end{proposition}

The proof will be presented in \cref{app:annealed-tight} for the sake of 
readability.

\subsection{Central Limit Theorem: Supercritical Case $\beta > 1$}\label{subsect:annealed-clt-supercrt}
The scaling we shall use in this supercritical case is different from the one 
we have considered so far. Envisioning the behavior of
\eqref{super-heuristic}, we introduce a new process $G^{\ep}_{\textsc{s}}$ 
defined by
\begin{equation}\label{def:G-ep-super-N}
G^\ep_{\textsc{s}} (t) := 
\ep^{\beta-\frac{1}{2}}
\Big(N^\ep(t) - \frac{\bar{m}_0(t)}{\ep^\beta}\Big).
\end{equation}
Then the equivalent of our decomposition \eqref{clt-decompose} is written as
\begin{equation}\label{eq:G-ep-Super-decomp}
G^\ep_{\textsc{s}}(t) =
\frac{\ep^\beta N^\ep(t) - m^\ep_0(t)}{\sqrt{\ep}}
+ \frac{m^\ep_0(t) -\bar{m}_0(t)}{\sqrt{\ep}}
:= G^\ep_{\textsc{s}, 1}(t)+ G^\ep_{\textsc{s}, 2}(t) \,,
\end{equation}
where $G^\ep_{\textsc{s}, 1}(t)$ and $G^\ep_{\textsc{s}, 2}(t)$ are respectively defined by (with $G^{\ep}_1 (t)$ introduced in \eqref{eqdef:Q-ep-1,2-fluctuation})
\begin{equation}\label{eq:G-ep-Super-decomp-exp}
G^\ep_{\textsc{s}, 1}(t) = \ep^{\frac{\beta - 1}{2}} G^{\ep}_1 (t), \quad \text{and} \quad
G^\ep_{\textsc{s}, 2}(t) =  \frac{m^\ep_0(t) -\bar{m}_0(t)}{\sqrt{\ep}} \,.
\end{equation}
As mentioned in \cref{subsect:intro-regimes-rescaledqueue}, the dominating term in our decomposition \eqref{eq:G-ep-Super-decomp} is $G^\ep_{\textsc{s}, 2}(t)$ whenever $\beta >1$ (this assertion will be established in the following propositions).
Furthermore, some of the results established in the subcritical case are still applicable here, so that we have
$G^\ep_{\textsc{s}, 1}(t)$ converging to $0$ in $\PP$-probability. Hence,
the analysis in supercritical case is completely captured by the
behavior of $G^\ep_{\textsc{s}, 2}(t)$. This global strategy will be detailed in the following propositions, starting from a univariate annealed CLT.

\begin{proposition}[Univariate Annealed CLT]\label{thm-FCLT-super-univeriate-annealed}
Assume \cref{hyp:arr-times-Gamma,hyp:F-bar-increments,hyp:clt-Z,supercri-hyp:coef-la-psi} are satisfied, and recall that $\beta > 1$ in this section.
Then for any fixed $t \geq 0$, the following annealed CLT holds true:
\begin{equation}\label{lim:G-ep-super-Gau-ann}
G^{\ep}_{\textsc{s}} (t) \Longrightarrow_{\PP} G_{\textsc{s}} (t)  \quad\quad \text{with} \quad
G_{\textsc{s}} (t) \overset{(\dd)}{=} \cn(0, \bar\sigma^2 (t)) \,,
\end{equation}
where the variance $\bar\sigma^2 (t)$ is defined by
\begin{equation}\label{def:barsig-sq-ann}
\bar\sigma^2 (t) = \iot \iot \la (s_1) \la (s_2) \bar F_{s_1} (t-s_1) \bar F_{s_2} (t-s_2) (s_1 \wedge s_2) \,\dd s_1 \dd s_2 \,.
\end{equation}
\end{proposition}

\begin{proof}
Recall the decomposition \eqref{eq:G-ep-Super-decomp} - \eqref{eq:G-ep-Super-decomp-exp} for $G^{\ep}_{\textsc{s}} (t)$. We will analyze the terms $G^\ep_{\textsc{s}, 1}(t)$ and $G^\ep_{\textsc{s}, 2}(t)$ therein separately. We divide the rest of the proof into several steps.

\noindent
{\it Step 1: $L^2$-limit for $G^\ep_{\textsc{s}, 1}(t)$.}
Recall from the expression in \eqref{eq:G-ep-Super-decomp-exp} that $G^\ep_{\textsc{s}, 1}(t) = \ep^{\frac{\beta - 1}{2}} G^{\ep}_1 (t)$ with $G^{\ep}_1 (t)$ introduced in \eqref{eqdef:Q-ep-1,2-fluctuation}. Moreover, we have $m^{\ep}_0 (t) = \ep^{\beta} m^{\ep} (t)$. Hence, one can write
$$
G^{\ep}_1 (t) = \ep^{\beta /2} \left(N^{\ep} (t) - m^{\ep} (t) \right) .
$$
Next we have seen that $N^{\ep} \sim \text{Poisson} (m^{\ep} (t))$. Therefore, conditioning on $Z$, it is readily checked that
\begin{equation}\label{def:G-ep-1-sq-var}
\E [|G^\ep_{1}(t)|^2]  =
\ep^{\beta} \E [|N^{\ep} (t) - m^{\ep} (t)|^2] = \ep^{\beta}  \text{Var} (N^{\ep} (t)) = \ep^{\beta} \Ebf \left[m^{\ep} (t)\right]
= \Ebf \left[m^{\ep}_0 (t) \right] \,.
\end{equation}
Next invoking \eqref{def:m-ep-expansionwith-m-ep-0} we have that 
\begin{equation}\label{def:m-ep0-with-h}
m^{\ep}_0 (t) = \iot h (s,t)\, \psi (Z_{s/\ep}) \,\dd s \,,
\quad \text{with} \quad h(s,t) = \la (s) \bar F_s (t-s) \,.
\end{equation}
Since $\la$, $\psi$ are bounded functions according to \cref{supercri-hyp:coef-la-psi} and $\int_{0}^{\infty} \bar F_s (u) \dd u$ is also bounded, it is easily seen that $m^{\ep}_0 (t) \leq C$ for a finite constant $C = C_{\la, \psi, \bar F}$, uniformly in $\ep$ and independently of the randomness in $Z$.
Plugging this information in~\eqref{def:G-ep-1-sq-var}, we have obtained
\begin{equation*}
\E [|G^\ep_{1}(t)|^2]  \leq C_{\la, \psi, \bar F} \,.
\end{equation*}
Hence recalling the relation in \eqref{eq:G-ep-Super-decomp-exp} for $G^\ep_{\textsc{s}, 1}(t) = \ep^{\frac{\beta - 1}{2}} G^{\ep}_1 (t)$ again and owing to the fact that $\beta >1$ in this supercritical case, we end up with
\begin{equation}\label{lim:G-eps1-L2}
\lim_{\ep \to 0} G^\ep_{\textsc{s}, 1}(t) = 0 \,, \quad\text{in} \quad L^{2} (\Omega, \PP) \,.
\end{equation}

\noindent
{\it Step 2: Convergence in distribution for $G^\ep_{\textsc{s}, 2}(t)$.}
Recall \eqref{eq:compute-hpsi-Z-int-var} for $m^{\ep}_0 (t) - \bar m_0 (t)$ with $h$ given in \eqref{def:m-ep0-with-h}. Hence, starting from the definition \eqref{eq:G-ep-Super-decomp-exp} for $G^\ep_{\textsc{s}, 2}(t)$, we get
\begin{equation*}
G^\ep_{\textsc{s}, 2}(t) =
h(t,t) 
\left(\ep^{-\frac{1}{2}} \iot \left(\psi(Z_{s/\ep}) -\bar{\psi} \right)\dd s
\right)
 - \int_0^t \frac{\dd}{\dd s}h(s,t) \left(
\ep^{-\frac{1}{2}} \int_0^s \left(\psi(Z_{r/\ep}) - \bar{\psi} \right)\dd r \right)
\dd s \,.
\end{equation*}
Let us now introduce a family of functionals defined on the space of continuous functions and indexed by $t \geq 0$:
\begin{equation}\label{def:functional-Phi-cont}
\phi_t (f) = h (t,t) f(t) - \iot \frac{\dd}{\dd s} h (s,t) f(s) \, \dd s \,.
\end{equation}
For any fixed $t \geq 0$, thanks to expression~\eqref{def:m-ep0-with-h} for $h$ and owing to the regularity properties of $\la$ and $\bar F$, it can be seen that $\phi_t: C(\R_+) \longrightarrow \R$ is continuous (recall that $C(\R_+)$ is equipped with the topology of uniform convergence on compact sets). Moreover, with this functional $\phi_t$ in hand, one can recast $G^\ep_{\textsc{s}, 2}(t)$ in the above expression as
\begin{equation}\label{eq:G-ep-S2-phi}
G^\ep_{\textsc{s}, 2}(t) = \phi_t (Y^{\ep}) \,, \quad \text{where}\quad
Y^{\ep} (s) =\frac{1}{\sqrt{\ep}} \ios \Big(\psi (Z_{r/\ep}) - \bar\psi \Big) \,\dd r \,.
\end{equation}
Therefore, invoking  \cref{hyp:clt-Z} we obtain that as $\ep \to 0$,
\begin{equation*}
G^\ep_{\textsc{s}, 2}(t) \Longrightarrow_{\PP} \phi_t (W) \, ,
\end{equation*}
where $W$ is a standard Brownian motion.

\noindent
{\it Step 3: Identification of the law.}
It is clear from expression \eqref{def:functional-Phi-cont} that $\phi_t (W)$ is a centered Gaussian random variable. Moreover, undoing the integration by parts~\eqref{eq:compute-hpsi-Z-int-var}, we have
\begin{equation*}
\phi_t (W) = \iot h (s,t) W (s) \,\dd s \,.
\end{equation*}
Since $\E [W (s_1) W(s_2)] = s_1 \wedge s_2$ for all $s_{1},s_{2}\ge 0$, we obtain the following expression for the quantity $\bar\sigma^2 (t)$ in~\eqref{lim:G-ep-super-Gau-ann}:
\begin{equation}\label{eq:var-G-eps-W}
\bar\sigma^2 (t)
:=
\E \left[ (\phi_t (W))^2\right] = \iot \iot h(s_1, t) h (s_2, t) (s_1 \wedge s_2) \, \dd s_1 \, \dd s_2 \,.
\end{equation}

\noindent
{\it Step 4: Conclusion.}
Going back to our decomposition \eqref{eq:G-ep-Super-decomp}-\eqref{eq:G-ep-Super-decomp-exp}, we have obtained
\begin{equation*}
G^\ep_{\textsc{s}, 1}(t) \xrightarrow{L^2 (\Omega)} 0 \quad \text{and} \quad
G^\ep_{\textsc{s}, 2}(t) \Longrightarrow_{\PP} \cn (0, \bar\sigma^2 (t)) \,.
\end{equation*}
Thus, relation \eqref{lim:G-ep-super-Gau-ann} follows, which completes the proof.
\end{proof}

With \cref{thm-FCLT-super-univeriate-annealed} in hand,
we may also describe the approximate self-similarity of the random variable $G_{\textsc{s}}$.
The approximate scaling for the random variable $G_{\textsc{s}}(t)$, similar to what we have obtained in \cref{prop:self-similar-process-fbm}, is given below. Notice the wider range of fluctuations in the supercritical regime.

\begin{proposition}\label{prop:self-similarity-G-s-bdd}
Let the assumptions of \cref{thm-FCLT-super-univeriate-annealed} prevail and consider the random variable $G_{\textsc{s}}(t)$ defined by \eqref{lim:G-ep-super-Gau-ann}, with variance $\bar\sigma^2 (t)$ spelled out in \eqref{def:barsig-sq-ann}. 
Then for $\al \in (0,1)$, there exist two constants $0 < c_1 < c_2 < \infty$ 
such that for all $t \geq 1$, we have
\begin{equation}\label{ineq:bounds-self-similarity-Gsann}
c_1 \, t^{3-2\al} \leq \bar\sigma^2 (t) \leq c_2 \, t^{3-2\al} \,,
\end{equation}
while for $\al > 1$, we have
\begin{equation}\label{ineq:bounds-self-similarity-Gsann-LT}
c_1 \, t \leq \bar\sigma^2 (t) \leq c_2 \, t.
\end{equation}
\end{proposition}

\begin{remark}
Notice that in the heavy tail scenario $\al \in (0,1)$, 
the exponent $3-2\al$ in \eqref{ineq:bounds-self-similarity-Gsann} is larger than $1$. This is in contrast with the exponent $1-\al$ obtained in 
\eqref{ineq:expected-G(t)} for the subcritical case. 
As mentioned in the introduction section, it indicates that 
in the heavy-traffic regime, the presence of random environment can induce 
some extra large fluctuations leading to superdiffusive 
behaviors.
\end{remark}

\begin{proof}
We start with the expression \eqref{def:barsig-sq-ann} for $\bar\sigma^2 (t)$. 
By assumption \eqref{def:F-bar-generalization} in \cref{hyp:F-bar-increments}, 
it can essentially be estimated as:
\begin{equation}\label{def:sigma-recast-selfsimilar}
\bar\sigma^2 (t) 
\asymp
\int_0^t\int_0^t
\left(\frac{1}{(t-s_1)^\al}\wedge 1\right)
\left(\frac{1}{(t-s_2)^\al}\wedge 1\right)s_1\wedge s_2
\, \dd s_1\, \dd s_2.
\end{equation}
The above integral can in fact be explicitly integrated. For this purpose,
we write
\begin{eqnarray*}
\bar\sigma^2 (t) 
& \asymp &
\int_{t-1}^t\int_{t-1}^t\,\cdots\,\,\dd s_2\,\dd s_1
+
\int_{t-1}^t\int_0^{t-1}\,\cdots\,\,\dd s_2\,\dd s_1
+
\int_0^{t-1}\int_0^{s_1}\,\cdots\,\,\dd s_2\,\dd s_1\\
& := &
I + II + III
\end{eqnarray*}
where the $\cdots$'s refer to the integrand in 
\eqref{def:sigma-recast-selfsimilar}
and without loss of generality, we have just considered ``half'' of the domain
of integraion: $0 < s_2 < s_1 < t$. Now we estimate $I, II$, and $III$.

For $I$, it is clear that $I\asymp t$.

For $II$, we compute,
\begin{eqnarray*}
II & = & 
\int_{t-1}^t\left[\int_0^{t-1}
\frac{s_2}{(t-s_2)^\al}\,\dd s_2\right]\,\dd s_1
= \int_0^{t-1} \frac{s_2}{(t-s_2)^\al}\,\dd s_2\\
& = &
\int_0^{t-1} 
t(t-s_2)^{-\al}
-(t-s_2)^{-\al+1}
\,\dd s_2\\
& = & 
\frac{t^{2-\al}}{(1-\al)(2-\al)}
-\frac{t}{1-\al}
+ \frac{1}{2-\alpha}.
\end{eqnarray*}

For $III$, we have
\begin{eqnarray*}
III & = & 
\int_0^{t-1}\int_0^{s_1} 
\frac{1}{(t-s_1)^\al} \frac{s_2}{(t-s_2)^\al} \,\dd s_2\,\dd s_1
= 
\int_0^{t-1}
\frac{1}{(t-s_1)^\al} 
\left[\int_0^{s_1} 
\frac{s_2}{(t-s_2)^\al} \dd s_2\right]
\,\dd s_1
\\
& = & 
\int_0^{t-1} 
\frac{1}{(t-s_1)^\al} 
\left[\int_0^{s_1} 
t(t-s_2)^{-\al}
-(t-s_2)^{-\al+1}
\,\dd s_2\right]
\,\dd s_1\\
\end{eqnarray*}
which after some tedious computation, gives
\[
III = 
\frac{t^{3-2\al}}{2(1-\al)^2(3-2\al)}
-
\frac{t^{2-\al}}{(1-\al)^2(2-\al)}
+\frac{t}{2(1-\al)^2}
-\frac{1}{(2-\al)(3-2\al)}.
\]

Estimates \eqref{ineq:bounds-self-similarity-Gsann} and \eqref{ineq:bounds-self-similarity-Gsann-LT} follow by summing up the above formulae
for $I$, $II$, and $III$.
\end{proof}

We are now ready to state and prove the functional version of our annealed CLT in the supercritical case. Its proof is similar to previous computations in this paper, and therefore will only be sketched.

\begin{theorem}[Functional CLT, Supercritical, Annealed]
\label{thm-FCLT-super-annealed}
Let $G^\ep_{\textsc{s}}$ be the process defined by \eqref{def:G-ep-super-N}, and assume \cref{hyp:arr-times-Gamma,hyp:F-bar-increments,hyp:clt-Z} are satisfied. Then we have
\begin{equation}\label{lim:Gs-convg-ann}
\left\{G^\ep_{\textsc{s}} \,;\, t \geq 0\right\} \Longrightarrow_{\PP} 
\left\{G_{\textsc{s}} \,;\, t \geq 0\right\} ,
\end{equation}
where $G_{\textsc{s}}$ is a continuous centered Gaussian process with covariance function
\begin{equation}\label{b1}
\E \left[G_{\textsc{s}} (s) G_{\textsc{s}} (t)\right] =
\iot \ios \la (r_1) \la (r_2) \bar F_{r_1} (s-r_1) \bar F_{r_2} (t-r_2) (r_1 \wedge r_2) \dd r_1 \dd r_2 \,.
\end{equation}
In \eqref{lim:Gs-convg-ann}, the convergence holds in the Skorokhod topology of c\`adl\`ag processes.
\end{theorem}

\begin{proof}
With \cref{thm-FCLT-super-univeriate-annealed} in hand, we are left with two tasks: (i) multivariate convergence; (ii) tightness in $D\left([0, \infty)\right)$. We will treat those two tasks sequentially below.

\noindent{\it Step 1: Multivariate Convergence.}
For the sake of conciseness, we will only deal with bivariate convergence and let the patient reader fill in the details for multivariate cases. Now consider two instants $s < t$ and the vector $(G^\ep_{\textsc{s}} (s), G^\ep_{\textsc{s}} (t))^T$. Thanks to \eqref{eq:G-ep-Super-decomp}, we obviously have
\begin{equation*}
(G^\ep_{\textsc{s}} (s), G^\ep_{\textsc{s}} (t))^T =
(G^\ep_{\textsc{s}, 1} (s), G^\ep_{\textsc{s}, 1} (t))^T + (G^\ep_{\textsc{s}, 2} (s), G^\ep_{\textsc{s}, 2} (t))^T \,.
\end{equation*}
Similarly to \eqref{lim:G-eps1-L2}, one can check that
\begin{equation*}
\lim_{\ep \to 0} (G^\ep_{\textsc{s}, 1} (s), G^\ep_{\textsc{s}, 1} (t)) = \mathbf{0}
 \,, \quad\text{in} \quad L^{2} (\Omega, \PP) \,.
\end{equation*}
Moreover, along the same lines as for \eqref{eq:G-ep-S2-phi}, we have
\begin{equation*}
(G^\ep_{\textsc{s}, 2} (s), G^\ep_{\textsc{s}, 2} (t)) = ( \phi_s (Y^{\ep}),  \phi_t (Y^{\ep}) )\,.
\end{equation*}
Since $(\phi_s, \phi_t): C(\R_+) \to \R^2$ is a continuous function, one can resort to \cref{hyp:clt-Z} again and obtain that as $\ep \to 0$,
\begin{equation*}
(G^\ep_{\textsc{s}, 2} (s), G^\ep_{\textsc{s}, 2}(t) ) \Longrightarrow_{\PP} (\phi_s (W), \phi_t (W)) \,.
\end{equation*}
Notice that the covariance function \eqref{b1} of $G_{\textsc{s}}$ is obtained by computing $\E \left[\phi_s (W) \phi_t (W)\right]$. This is achieved in the same way as for \eqref{eq:var-G-eps-W}.

\noindent{\it Step 2: Tightness.}
Here, we make use of decomposition \eqref{eq:G-ep-Super-decomp} again. Note that in the subcritical case, $G^{\ep}_1 (t)$ is the dominating term, as illustrated heuristically in \eqref{sub-heuristic}. 
Now in the current supercritical case, the same argument as in 
\cref{prop:tightness-annclt} shows that the family of processes $\{G^\ep_{\textsc{s}, 1} \,;\, \ep>0\}$ is tight in $D ([0, \infty))$.

As far as $G^\ep_{\textsc{s}, 2}$ is concerned, we invoke the relation $G^\ep_{\textsc{s}, 2}(t) = \phi_t (Y^{\ep})$ in \eqref{eq:G-ep-S2-phi} and then define a global application
\begin{equation*}
\phi: C(\R_+) \longrightarrow C(\R_+) \,, \quad \text{with} \quad \phi (f) = \left\{\phi_t (f) \,;\, t \geq 0 \right\} .
\end{equation*}
In addition, owing to \eqref{def:functional-Phi-cont} and thanks to the fact that $h$ is a bounded function, we have
\begin{equation}\label{phi-bound-sup}
\|\phi(f)\|_{L^\infty [0,T]} \leq C_{\la, \bar F, T} \, \|f\|_{L^\infty [0,T]} \,.
\end{equation}
The upper bound \eqref{phi-bound-sup} together with the fact that $\phi$ is
a continuous map gives tightness in $C([0, T])$ for all $T>0$,
which in turn implies the tightness of $\{G^\ep_{\textsc{s}, 2} \,;\, \ep>0\}$ in the space $D([0, \infty))$.

Combining the tightness of $G^\ep_{\textsc{s}, 1}$ and $G^\ep_{\textsc{s}, 2}$, we have thus obtained the tightness of the whole sequence $G^\ep_{\textsc{s}}$, which completes the proof.
\end{proof}

The continuity result for the process $G_{\textsc{s}}$ mirrors 
\cref{G-continuity}.

\begin{proposition}\label{superG-continuity}
Let $G_{\textsc{s}}$ be the centered Gaussian process given in \eqref{lim:Gs-convg-ann} with covariance function \eqref{b1}.
Then $G_{\textsc{s}}$ is a $\mu$-H\"older continuous process on any interval $[0,T]$, for all $\mu \in (0, 1/2)$.
\end{proposition}

\begin{proof}
Since $G_{\textsc{s}}$ is a Gaussian process, it suffices to show that for $s < t$,
\begin{equation*}
\E \left[(G_{\textsc{s}} (t) - G_{\textsc{s}} (s))^2\right] \,\lesssim\, t-s \,.
\end{equation*}
Then the proof of this claim follows along the same lines as that of \cref{G-continuity}.
\end{proof}

\begin{remark}
As mentioned in \cref{subsect:intro-regimes-rescaledqueue}, the limiting behavior of $G^\ep$ that is defined by \eqref{centered-Nep} is determined by the dominating term. Namely, we have analyzed $G^\ep = G^\ep_1 + G^\ep_2$ as shown in \eqref{clt-decompose}-\eqref{eqdef:Q-ep-1,2-fluctuation}. In the subcritical case $\beta < 2 \ka$, the dominating term is $G^\ep_1$. In contrast, the dominating term becomes $G^\ep_2$ in the supercritical case $\beta > 2\ka$. The critical situation $\beta = 2\ka$ is more intricate in that the two processes $G^\ep_1$ and $G^\ep_2$ are of equal magnitude. In addition, they are far from being independent. For the sake of conciseness, we have thus deferred the analysis of this delicate situation to a subsequent publication. See also the discussion in the
Introduction.
\end{remark}

\section{Summary and Future Work.}
We considered infinite-server queues driven by $\cox$ (or doubly stochastic Poisson) process input in a fast oscillatory (or rapidly averaging) random environment.  Exact analysis of these queueing models poses difficulties, and therefore approximations are warranted. We establish diffusion approximations to the (re-scaled) number-in-system process in both the quenched and annealed settings, by establishing FCLTs using a unified stochastic homogenenization approach. This contrasts and complements the heavy-traffic frameworks that have been used extensively to prove FCLTs in the annealed setting (alone). Furthermore, our analysis uncovers  sub- and supercritical regimes in which the annealed diffusion approximations are qualitatively different (even though the arrival and service primitives both are identical in the two regimes). Furthermore, we show that quenched FCLTs cannot be proved in the supercritical regime. 

Left out of the existing analysis is the so-called critical regime where the analysis is more delicate. We leave this to a subsequent paper. The qualitative differences in the annealed diffusion approximations alluded to above are potentially significant for system performance  optimization and analysis. This too is left to a future paper to be studied in depth. Finally, our analysis can be used to analyze other queueing models and networks, such as generalized Jackson networks with Cox process input.

\appendix
\section{Proof of \cref{prop:CLT-G-ep-md}: \\ Multivariate CLT in the
Quenched Setting}\label{sect:appendix-prop-multclt}

\begin{proof}[Proof of \cref{prop:CLT-G-ep-md}]
The proof follows the same lines as its bivariate counterpart. 
For the sake of clarity, in this proof,
we use the convention that
$(\cdot)_{i,j}$ refers to the dependence on $(t_i, t_j$) with $t_i < t_j$.
The proof is divided into several steps.

\noindent
{\it Step 1: Disjoint Regions.}
In the same spirit as the bivariate case, we shall decompose the quadrant $\R_+ \times \R_+$ into several disjoint regions: 
\[
\{A_{i, j} : \text{for} \,\, 1\leq i \leq n \,\, \text{and}\,\, 1 \leq j \leq n-i+1 \} \,,
\]
where for each $i$ and $j$, the region is defined as
\begin{equation}\label{defeq:A-ij-disjointregions}
A_{i, j} := \left\{(\gamma, l) \in \R_+ \times \R_+: \, t_{i-1} \leq \gamma \leq t_i \,\,\text{and} \,\,t_j \leq \gamma + l \leq t_{j+1} \right\} ,
\end{equation}
with the additional convention $t_0 = 0$ and $t_{n+1} = \infty$. Since the $A_{i,j}$’s are disjoint regions, the quantities $\{M^{\ep} (A_{i,j}) : i, j = 1,2, ... , n\}$ are independent Poisson random variables. Similar to \eqref{eq:m-t-def}-\eqref{def:m-ep-expansionwith-m-ep-0} and \eqref{eq:def-mean-M-ep-A_i-reg}, their respective quenched means are given by
\begin{equation*}
m^{\ep}(A_{i,j}) = \E_Z [M^{\ep} (A_{i,j})] 
= \frac{1}{\ep^{\beta}} \int_{A_{i,j}} \nu (s, \dd r) \la (s) \psi (Z_{s/\ep}) \dd s \,.
\end{equation*}

\noindent
{\it Step 2: CLT for $M^{\ep} (A_{i, j})$.}
A central limit theorem for each family $\{M^{\ep} (A_{i,j}) \,;\, \ep >0\}$ can be obtained along the same lines as in \eqref{defeq:M-ep-A-i-m-ep-i}-\eqref{eq:note-for-recaled-M-ep}. Indeed, according to \cref{prop:gen-limit-m-ep} case (ii), we know that $\Pbf$-almost surely, 
$\lim_{\ep \to 0} m^{\ep}(A_{i,j}) = \infty$. In other words, it is easily seen that
\begin{equation*}
\lim_{\ep \to 0} \ep^{\beta} m^{\ep}(A_{i,j})
= \lim_{\ep \to 0} \ep^{\beta} \E_Z [M^{\ep} (A_{i,j})] =\bar m_0(A_{i,j}) \,.
\end{equation*}
Recalling (see \cref{hyp:F-bar-increments}) that $\nu$ admits a density $\ell$, we obtain
\begin{equation*}
\bar m_0(A_{i,j}) = \bar\psi \int_{A_{i,j}} \la (s) \, \ell_s (r) \, \dd r \,\dd s
\,.
\end{equation*}
By following the proof as in the bivariate case,
we obtain the following central limit theorem for $M_{i,j}^{\ep}$:
\begin{equation}\label{eq:clt-M-i-j-scaled}
\widetilde{M^\ep}(A_{i,j}) \equiv \ep^{\beta /2} 
\left(M^{\ep} (A_{i,j}) - m^{\ep}(A_{i,j}) \right) 
\xrightarrow{(\dd)} \mathcal{N} (0, \bar m_0(A_{i,j})).
\end{equation}


\noindent
{\it Step 3: Multidimensional CLT.}
We now recall that $G^{\ep} (t)$ is defined by \eqref{centered-Nep}, as well as the limit shown in \eqref{eq:clt-M-i-j-scaled}. Due to the definition in \eqref{defeq:A-ij-disjointregions} of the regions $A_{i,j}$, it can be easily seen that for all $k$ such that $1 \leq k \leq n$,
\begin{equation}\label{eq:G-ep-express-in-M_ij}
G^{\ep} (t_k) = \sum_{i = 1}^{k} \sum_{j = k-i+1}^{n-i+1} 
\widetilde{M^\ep}(A_{i,j}) \,.
\end{equation}
Due to the fact that the random variables $M_{i,j}^{\ep}$ defined by \eqref{eq:clt-M-i-j-scaled} are independent, we will express $G^{\ep} (t_k)$’s in terms of $M^\ep(A_{i,j})$’s.
Specifically, we note that the total number of the regions denoted by $A_{i,j}$ is $n(n+1)/2$.
Then according to \eqref{eq:G-ep-express-in-M_ij}, we can decompose the vector $\mathbf{u} \equiv (G^{\ep} (t_1), G^{\ep} (t_2), \ldots , G^{\ep} (t_n))^T$ as follows
\begin{multline}\label{eq:matrixmult-uvQmix}
\mathbf{u} = Q_{mix} \mathbf{v}\,, \quad \text{with}\\
\mathbf{v} = \left(\widetilde{M^\ep}(A_{1,1}) , \widetilde{M^\ep}(A_{1,2}) , 
\ldots , \widetilde{M^\ep}(A_{1,n}), 
\widetilde{M^\ep}(A_{2,1}), \ldots  \widetilde{M^\ep}(A_{2,n-1}), \ldots , 
\widetilde{M^\ep}(A_{n,1})\right)^T ,
\end{multline}
where $\mathbf{u} \in \R^{n \times 1}$, $\mathbf{v} \in \R^{(n(n+1)/2) \times 1}$, while $Q_{mix} \in \R^{n \times (n(n+1)/2)}$ is a mixing matrix that only contains $0$ and $1$ (and whose exact expression is omitted for the sake of conciseness).
Moreover, each entry in the vector $\mathbf{v}$ defined by \eqref{eq:matrixmult-uvQmix} satisfies the CLT \eqref{eq:clt-M-i-j-scaled}. Those entries are independent, owing to the fact that the regions $A_{i,j}$ are disjoint. Therefore, we have obtained that as $\ep \to 0$,
\begin{equation}\label{lim:matrix-v-M-ep-with-Imbar}
\mathbf{v} \xrightarrow{(\dd)}
\mathcal{N} \left(0, I_{\bar m} \right) \quad \text{with} \quad I_{\bar m} =
\text{Diag} \left(
\bar m_0(A_{1,1}) , \bar m_0(A_{1,2}) , \ldots , \bar m_0(A_{n,1}) \right),
\end{equation}
namely, $I_{\bar m} \in \R^{(n(n+1)/2) \times (n(n+1)/2)}$ contains the diagonal entries $\bar m_0(A_{i,j})$ and $0$ everywhere else. 
We now invoke some standard considerations on Gaussian vectors in order to conclude that the vector $\mathbf{u}$ defined by \eqref{eq:matrixmult-uvQmix} satisfies the following limit in distribution with respect to $\PP_Z$:
\begin{equation*}
\mathbf{u}
\xrightarrow{(\dd)}
\mathcal{N} (0, \Gamma) \,,
\quad \text{where} \quad
\Gamma = Q_{mix} \,\, I_{\bar m} \,\, Q_{mix}^{T} \in \R^{n \times n}\,.
\end{equation*}

One could easily get expression \eqref{defeq:Gamma-titj-inmultdim} for the matrix $\Gamma$ starting from \eqref{lim:matrix-v-M-ep-with-Imbar}. However, it might be shorter to argue as follows: the $n$-dimensional limiting distributions given by \eqref{lim:mult-G-ep-t1-tn} have to be compatible with the $2$-dimensional limiting distributions in \eqref{eq:clt-for-G-ep-vec}. Therefore according to \cref{rmk:expression-Gamma-titj}, each entry of the matrix $\Gamma$ in \eqref{lim:mult-G-ep-t1-tn} is given by the expression~\eqref{eq:gamma-titj-with-psi-bar} in \cref{rmk:expression-Gamma-titj}. This proves \eqref{defeq:Gamma-titj-inmultdim}.
\end{proof}

\section{Proof of \cref{prop:G-ep-family-tightness}:\\
Tightness in the Quenched Setting}~\label{app:quenched-tight}
As in \cite{RR}, we shall prove the tightness of our sequence of processes $\{G^{\ep} \,;\, \ep >0\}$ in the space of c\`adl\`ag functions $D\left([0, \infty)\right)$. Our strategy will be based on considerations concerning this functional space, as well as Orlicz-type norms for random variables. Let us start by giving a proper definition of Orlicz norms.

\begin{definition}\label{defnthm:orlicz-norm}
Let $X$ be a random variable defined on a probability space $(\Omega, \mathcal F, \PP)$. Also consider the function $\phi$ defined by
\begin{equation}\label{def:phi-fun}
\phi (x) = e^x -1, \quad x \geq 0 \,.
\end{equation}
Then the Orlicz norm $\norm{X}$ of $X$ is given by
\begin{equation}\label{def:orlicz-norm}
\norm{X} = \inf \left\{c > 0: \E \lc\phi (c^{-1} |X|)\rc \leq 1 \right\} .
\end{equation}
Notice that $\norm{X}$ can be rewritten as
\begin{equation}\label{def:orlicz-norm-ver}
\norm{X} = \inf \left\{c>0 : \E \left[e^{\frac{|X|}{c}} \right] \leq 2 \right\} .
\end{equation}
\end{definition}

For our quenched analysis, we shall need some quenched Orlicz norms. We thus label the notation below for further use.

\begin{notation}\label{notation:orlicz-norm-Z}
We denote by $\norm{X}_Z$ the Orlicz norm \eqref{def:orlicz-norm} when the quenched probability $\PP_Z$ is considered. Namely, we set
\begin{equation*}
\norm{X}_Z = \inf \left\{c > 0: \E_Z \lc\phi (c^{-1} |X|)\rc \leq 1 \right\}
= \inf \left\{c>0 : \E_Z \left[e^{\frac{|X|}{c}} \right] \leq 2 \right\} .
\end{equation*}
\end{notation}

With \cref{defnthm:orlicz-norm,notation:orlicz-norm-Z} in hand, we state a lemma that is an adaptation of \cite[Lemma 1]{RR} to the quenched setting.

\begin{lemma}\label{lemma:Poisson-rv-Orlicz}
Let $N^{\ep} (t)$ be the conditional Poisson random variable defined by \eqref{eq:N-t1}, with quenched mean $m^{\ep} (t)$ for a given strictly positive $\beta$. Then there exists a set $\Omega$ with $\Pbf (\Omega) = 1$ satisfying the following: for all $\omega \in \Omega$, there is an $\ep_0 (\omega) >0$ such that  for all $\ep \leq \ep_0$ we have a deterministic constant $K >0$ such that
\begin{equation*}
\mathbbm{1}_{\{m^{\ep} (t) \geq 1\}}
\norm{N^{\ep} (t) - m^{\ep} (t)}_Z \leq K \sqrt{m^{\ep} (t)} \,.
\end{equation*}

\end{lemma}

\begin{proof}
As mentioned above, this is a mere quenched version of \cite[Lem\-ma~1]{RR}.
To apply this result, first, we check the existence of the set $\Omega$ with $\Pbf (\Omega) = 1$ such that for all $\omega \in \Omega$, there is an $\ep_0 (\omega) >0$ and for all $\ep \leq \ep_0$ we have $m^{\ep} (t) \geq 1$. 
Indeed, this follows from~\eqref{eq:limit-rescaled-m-ep} which asserts that 
for $\beta >0$, we have $\lim_{\ep \to 0} m^{\ep} (t) = \infty$. 
Second, 
note that even though $m^\ep(t)$ is a random variable depending on 
a realization of $Z$, the desired statement is in essence a property of
Poisson random variable.

We now compute the constant bound $K$.
With $m^{\ep} (t) \geq 1$, the inequality below for any $Z = Z(\omega)$ immediately follows from \cite[Lem\-ma~1]{RR}:
\begin{equation*}
\norm{N^{\ep} (t) - m^{\ep} (t)}_Z \leq K_Z \sqrt{m^{\ep} (t)} \,,
\end{equation*}
where the random variable $K_Z$ can be written as $K_Z = K_{Z,1} + K_{Z,2}$ with
\begin{eqnarray}
K_{Z,1} &=& \frac{1}{\sqrt{2 e^2 \log (1+\PP_Z [N^{\ep} (t) > m^{\ep} (t)])}} \,, \label{eq:const-orlicz-Poisson1}\\
K_{Z,2} &=& \frac{1}{\sqrt{2 e^2 \log (1+\PP_Z [N^{\ep} (t) \leq m^{\ep} (t)])}} \,. \label{eq:const-orlicz-Poisson2}
\end{eqnarray}
In order to bound the random variable $K_Z$, let $\Pcal$ denote any generic probability measure. Then we consider $p_1$ and $p_2$ defined by
\begin{equation*}
p_1 = \sup_{m \geq 1} \Pcal [\text{Poisson} (m) \leq m] \in (0,1) 
\quad \text{and} \quad
p_2 = \sup_{m \geq 1} \Pcal [\text{Poisson} (m) > m)] \in (0,1) \,.
\end{equation*}
Then we introduce some constants $K_1$ and $K_2$ defined as
\begin{equation*}
K_1 = \frac{1}{\sqrt{2 e^2 \log (1+ (1-p_1))}} \quad \text{and} \quad
K_2 = \frac{1}{\sqrt{2 e^2 \log (1+ (1-p_2))}} \,\,.
\end{equation*}
It is easily seen that $K_{Z,1} \leq K_1$ and $K_{Z,2}\leq K_2$. We now set $K=K_1 + K_2$ and we get the inequality of \cref{lemma:Poisson-rv-Orlicz}.
This completes our proof.
\end{proof}

The proof of the tightness involves an estimation of the quantities $\ep^{\beta} m^{\ep}_1$ and $\ep^{\beta} m^{\ep}_3$ in~\eqref{ineq:g-ep-j-bound} below.
Recall the expressions for $m^{\ep}_1, m^{\ep}_3$ are given in 
\eqref{eq:def-mean-M-ep-A_i-reg}. Since we are now interested in quantifying increments of the form $G^{\ep} (t_2) - G^{\ep} (t_1)$, to emphasize the
dependence on $t_1$ and $t_2$, 
we will set $m^{\ep}_i = m^{\ep}_i (t_1, t_2)$ and $A_i = A_i (t_1,t_2)$ in the definitions \eqref{def:region_A1}-\eqref{def:region_A3} involved in \eqref{eq:def-mean-M-ep-A_i-reg}. Therefore, taking into account the expression $\nu (s, \dd r) = \ell_s (r) \dd r$ given in \cref{hyp:F-bar-increments}, one can recast \eqref{eq:def-mean-M-ep-A_i-reg} for $i = 1,3$ as
\begin{equation}\label{eq:ep-m-1,3-int}
\ep^{\beta} m^{\ep}_i (t_1, t_2) = \int_{A_i (t_1,t_2)} \nu (s, \dd r) \la (s) \psi (Z_{s/\ep})\, \dd s
= \int_{A_i(t_1,t_2)} \ell_s (r) \la (s) \psi (Z_{s/\ep}) \,\dd r \,\dd s \,.
\end{equation}
See \cref{fig:A_i-region-add} for an illustration.
We state a lemma regarding uniform bounds of the quantities $\ep^{\beta} m^{\ep}_1 (t_1, t_2)$ and $\ep^{\beta} m^{\ep}_3 (t_1, t_2)$ for further use.

\begin{lemma}\label{lemma:unifbound-epbeta-m-i}
Let \cref{hyp:ergodicity,hyp:arr-times-Gamma,hyp:F-bar-increments} prevail.
In addition, we set
\begin{equation*}
\|\psi\|_{L^\infty(\R^d)} \equiv K < \infty.
\end{equation*}
Then there exists a constant $C$ such that for $0 \leq t_1 \leq t_2 \leq 1$, the quantities $\ep^{\beta} m^{\ep}_i (t_1, t_2)$ defined by
\eqref{eq:ep-m-1,3-int} satisfy
\begin{equation}\label{ineq:epbeta-mep-bounds}
\ep^{\beta} m^{\ep}_1 (t_1, t_2)+ \ep^{\beta} m^{\ep}_3 (t_1, t_2) \leq C \left|t_2-t_1 \right| . 
\end{equation}
\end{lemma}

\begin{proof}
Let us first handle the term $\ep^{\beta} m^{\ep}_1(t_1, t_2)$ in \eqref{eq:ep-m-1,3-int}. Owing to the bound \eqref{def:F-bar-generalization} on $\ell$ and to the fact that $\la$ is bounded on the interval $[0,1]$, there exists a constant $C >0$ (which can change from line to line) such that
\begin{equation*}
\ep^{\beta} m^{\ep}_1 (t_1, t_2) \leq C K \int_{A_1 (t_1, t_2)} \left( \frac{1}{r^{1+\al}} \wedge 1\right) \dd r \, \dd s \,.
\end{equation*}
On $A_1 (t_1, t_2)$ we will simply bound $(r^{-(1+\al)} \wedge 1)$ by $1$. 
This yields
\begin{equation}\label{ineq:bdd-epbetam1}
\ep^{\beta} m^{\ep}_1 (t_1, t_2) \leq C \,\text{Area} (A_1(t_1, t_2)) 
\leq C t_1(t_2-t_1) 
\leq C (t_2-t_1) 
\,,
\end{equation}
where the second inequality is easily seen from the definition \eqref{def:region_A1} of $A_1(t_1, t_2)$.

The term $\ep^{\beta} m^{\ep}_3(t_1, t_2)$ is treated similarly. Namely from \eqref{eq:ep-m-1,3-int} and the aforementioned bounds on $\la$, $\psi$, and $\ell$, we have
\begin{equation*}
\ep^{\beta} m^{\ep}_3 (t_1, t_2)\leq C K \int_{A_3(t_1, t_2)} \left( \frac{1}{r^{1+\al}} \wedge 1\right) \dd r \, \dd s \,.
\end{equation*}
Thanks to the definition \eqref{def:region_A3} of $A_3(t_1, t_2)$, we have
\begin{equation*}
\ep^{\beta} m^{\ep}_3 (t_1, t_2)\leq C \int_{t_1}^{t_2} \int_{0}^{\infty} \left( \frac{1}{r^{1+\al}} \wedge 1\right) \dd r \, \dd s \,,
\end{equation*}
and since $\int_{0}^{\infty} \left( r^{-(1+\al)} \wedge 1\right) \dd r $ is finite, we end up with
\begin{equation}\label{ineq:bdd-epbetam3}
\ep^{\beta} m^{\ep}_3 (t_1, t_2)
\leq C (t_2 - t_1) \,.
\end{equation}

Putting together \eqref{ineq:bdd-epbetam1} and \eqref{ineq:bdd-epbetam3}, our claim \eqref{ineq:epbeta-mep-bounds} is proved.
\end{proof}

We now label a lemma adapted from \cite[Section 3.1 Lemma 2]{RR}. We continue to denote by $K$ (sometimes with a subscript) any generic constant. The specific value of $K$ may change from line to line. 
Since we will deal with time points from the discrete set 
\eqref{def:F_j-dyadic}, for simplicity, with
$t_1 = 2^{-j} (k-1)$ and $t_2 = 2^{-j} k$ (referring to \cref{fig:A_i-region-add}) we write the dependence on 
$(t_1,t_2)$ as $[j,k]$. For example,
\[
A_i (t_1, t_2) = A_i [j,k]
\quad\text{and}\quad m^{\ep}_i (t_1, t_2) = m^{\ep}_i [j,k].
\]

\begin{lemma}\label{lemma:m-j-Orlicz-norm-bounds}
For $\ep >0$, let $G^{\ep}$ be defined as in \cref{prop:G-ep-family-tightness}. For $j \geq 0$, we set
\begin{equation}\label{def:m-j-max-G-ep-diff}
g^{\ep}_j = \bigvee_{k=1}^{2^j} \left|G^{\ep} (2^{-j} k) - G^{\ep} (2^{-j} (k-1))\right| .
\end{equation}
As in \cref{lemma:Poisson-rv-Orlicz}, we assume $\beta >0$.
Then there exists an $\ep_0 (\omega, j)$ such that for all $\ep \leq \ep_0 (\omega, j)$ we have
\begin{equation}\label{ineq:g-ep-j-bound}
\norm{g^{\ep}_j}_Z \leq K \sup_{1 \leq k \leq 2^j }\left(\sqrt{\ep^{\beta} m^{\ep}_1[j,k]} + \sqrt{\ep^{\beta} m^{\ep}_3[j,k]} \right) j \,,
\end{equation}
where $K$ is a constant defined in \cref{lemma:Poisson-rv-Orlicz}.
\end{lemma}

\begin{proof}
This proof is inspired by \cite[Section 3.1 Lemma 2]{RR}. For the sake of clarity and the reader's convenience, we present a complete proof as follows.

Referring to \cref{fig:A_i-region-add} for an illustration, let us set $t_1 = 2^{-j} (k-1)$ and $t_2 = 2^{-j} k$ for a given $k$ in $\{1, \ldots, 2^j\}$.
Then with the notation of \cref{prop:bivariate-clt-for-G-ep} and its proof, 
it is readily checked from \eqref{eq:G-ep-decompmatrix-M-ep} that
\begin{equation}\label{conc.G.diff}
G^{\ep} (2^{-j} k) - G^{\ep} (2^{-j} (k-1)) 
= 
\widetilde{M^{\ep}}(A_3[j,k]) -\widetilde{M^{\ep}}(A_1[j,k]) 
\end{equation}
where the random variables $\widetilde{M^{\ep}} (A_i)$ are introduced in \eqref{defeq:M-ep-A-i-m-ep-i}.
Hence with the Orlicz norm introduced in \cref{defnthm:orlicz-norm,notation:orlicz-norm-Z}, we get
\begin{multline*}
\norm{G^{\ep} (2^{-j} k) - G^{\ep} (2^{-j} (k-1)) }_Z \\
\leq \ep^{\beta /2} \norm{M^{\ep}(A_1[j,k]) - m^{\ep}_1[j,k]}_Z + 
\ep^{\beta /2} \norm{M^{\ep}(A_3[j,k]) - m^{\ep}_3[j,k]}_Z \,.
\end{multline*}
Then we estimate the terms $\norm{M^{\ep}(A_i [j,k]) - m^{\ep}_i [j,k]}_Z$ for $i=1,3$. To this aim, we notice that $\ep_0 = \ep_0 (\omega, j, k)$ is chosen to be small enough, and hence $m^{\ep}_i [j,k]\geq 1$ for all $\ep \leq \ep_0$.
Therefore, one can apply \cref{lemma:Poisson-rv-Orlicz} to get the following bound:
\begin{equation}\label{ineq:bounds-Orlicz-G-ep-increment-with-k}
\norm{G^{\ep} (2^{-j} k) - G^{\ep} (2^{-j} (k-1)) }_Z \leq
K \left(\sqrt{\ep^{\beta} m^{\ep}_1[j,k]} + \sqrt{\ep^{\beta} m^{\ep}_3[j,k]} \right) . 
\end{equation}

To achieve \eqref{ineq:g-ep-j-bound}, 
we resort to the fact that Orlicz norms are tailored to be well behaved under sup-type manipulations. Recall that we are considering the function $\phi$ defined by~\eqref{def:phi-fun}.
Then under the quenched probability $\PP_Z$, \cite[Lemma 2.2.2 p.96]{VW} asserts that
when $\ep_0$ is small enough, we have
\begin{equation}\label{ineq:bounds-Orlicz-G-ep-increment-with-k2}
\norm{ \bigvee_{k=1}^{2^j}  \left| G^{\ep} (2^{-j} k) - G^{\ep} (2^{-j} (k-1)) \right| }_Z
\leq C_{\phi} \, \phi^{-1} (2^j) \bigvee_{k=1}^{2^j} \norm{ G^{\ep} (2^{-j} k) - G^{\ep} (2^{-j} (k-1)) }_Z \,,
\end{equation}
for a constant $C_{\phi}$ depending only on $\phi$. 
In addition, it is easily seen that $\phi^{-1} (x) = \ln (x+1)$ for $x \geq 0$. Hence, plugging \eqref{ineq:bounds-Orlicz-G-ep-increment-with-k} into the above, our claim \eqref{ineq:g-ep-j-bound} is proved.
\end{proof}

We now turn to the proofs of the two main points in order to reach tightness, 
that is, inequality~\eqref{ineq:prob-sup-diff-G-ep} about increments of $G^\ep$ on the grid $F_v$ and~\eqref{lim:max-increment-G-ep-v} about the small scale
increments within a grid.

\begin{lemma}[Proof of \eqref{ineq:prob-sup-diff-G-ep}]\label{lemma:tightness-ineqstate1}
As in \cref{prop:G-ep-family-tightness}, we consider the family $\{G^{\ep} \,;\, \ep >0\}$  of processes such that $G^{\ep} (t)$ is defined by \eqref{centered-Nep}.
Then for all rational $\theta ,  \eta >0$, one can find a set $\Omega_{\theta, \eta}$ with $\Pbf (\Omega_{\theta, \eta}) = 1$ and for all $\omega \in \Omega_{\theta, \eta}$, there exist $\ep_0 (\omega),  \delta (\omega) >0$
such that~\eqref{ineq:prob-sup-diff-G-ep} holds true: namely for all $\ep \leq \ep_0 (\omega)$, there exists $v_0 = v_0 (\ep_0 (\omega))$ such that for $v \geq v_0$ we have
\begin{equation*}
\PP_Z \left[ \sup_{t_1, t_2 \in F_v \,;\, |t_2-t_1| \leq \delta} |G^{\ep} (t_2) - G^{\ep} (t_1)| > \theta\right]
< \eta \,,
\end{equation*}
where $F_v$ is defined by \eqref{def:F_j-dyadic}.

\end{lemma}

\begin{proof}
This proof is inspired by \cite{RR} but will be adapted to our quenched settings under $\Pbf$-almost sure statement.
For all rational $\theta ,  \eta >0$ and $\omega \in \Omega_{\theta, \eta}$ with $\Pbf (\Omega_{\theta, \eta}) = 1$, we are allowed to choose $\ep \leq \ep_0 (\omega)$ and $v \geq v_0$ such that \eqref{ineq:prob-sup-diff-G-ep} holds true.
Recall from \eqref{def:anneal-P} that $\Pbf$ denotes the probability distribution of $Z = Z(\omega)$.
In the sequel we will define $v = v(\ep)$ such that as $\ep \to 0$, $v \to \infty$.

Fix $i$ and assume $\ep_0 (\omega)$ is small enough to make $v \geq i$ for all $\ep \leq \ep_0$. For $j \geq 1$, let us recall the definition \eqref{def:F_j-dyadic} for $F_j$. For $t \in [0,1]$ and $i \leq j \leq v$, we define $t_j := \sup \{k 2^{-j} \in F_j: k 2^{-j} \leq t\}$.
Hence, $t \in F_v$ implies that $t_{j-1} = t_j$ or $t_{j-1} = t_j - 2^{-j}$. Also notice that $|t-s| \leq 2^{-i}$ implies that $t_i = s_i$ or $t_i = s_i \pm 2^{-i}$.
Then we write
\begin{equation}\label{eq:G-ep-t-expansionform}
G^{\ep} (t) = \sum_{j=i+1}^{v} (G^{\ep} (t_j) - G^{\ep} (t_{j-1})) + G^{\ep} (t_i) \,,
\end{equation}
with a similar expression for $G^{\ep} (s)$. Now consider $s, t \in F_v$ with $|t-s| \leq 2^{-i}$. Owing to \eqref{eq:G-ep-t-expansionform}, we have
\begin{equation*}
|G^{\ep} (t) - G^{\ep} (s)| \leq \sum_{j=i+1}^{v} |G^{\ep} (t_j) - G^{\ep} (t_{j-1})| + \sum_{j=i+1}^{v} |G^{\ep} (s_j) - G^{\ep} (s_{j-1})| +|G^{\ep} (t_i) - G^{\ep} (s_i)| \,.
\end{equation*}
Therefore, with the definition \eqref{def:m-j-max-G-ep-diff} of $g^{\ep}_j$ in mind, it is easily seen that
\begin{equation*}
|G^{\ep} (t) - G^{\ep} (s)|
\leq \left(2 \sum_{j=i+1}^{v} g^{\ep}_j \right) +g^{\ep}_i \,\, \leq \,\, 2 \sum_{j=i}^{v} g^{\ep}_j \,.
\end{equation*}
This bound does not depend on $s, t \in F_v$ with $|t-s| \leq 2^{-i}$. Hence, we have obtained
\begin{equation}\label{eqdef:xi-ep-i}
\xi_{\ep , i} := \bigvee_{s, t \in F_v , |t-s| \leq 2^{-i}} |G^{\ep} (t) - G^{\ep} (s)| \leq 2 \sum_{j=i}^{v} g^{\ep}_j \,.
\end{equation}
Recall that in \cref{lemma:m-j-Orlicz-norm-bounds}, we have bounded $\norm{g^{\ep}_j}_Z$ for all $i \leq j \leq v$ by considering a sufficiently small parameter $\ep \leq \ep_0 (\omega)$. Therefore, we obtain
\begin{eqnarray}
\norm{\xi_{\ep , i}}_Z &=& \norm{\bigvee_{s, t \in F_v , |t-s| \leq 2^{-i}} |G^{\ep} (t) - G^{\ep} (s)|}_Z \leq 2 \sum_{j=i}^{v} \norm{g^{\ep}_j}_Z \notag \\
&\leq& 2K \sum_{j=i}^{v}  j \sup_{1\leq k\leq2^j} \left(\sqrt{\ep^{\beta} m^{\ep}_1[j,k]} + \sqrt{\ep^{\beta} m^{\ep}_3[j,k]} \right) .\label{ineq:xi-ep-i-norm}
\end{eqnarray}
Let us apply \cref{lemma:unifbound-epbeta-m-i} to the right hand side of \eqref{ineq:xi-ep-i-norm}. We recall that $m^{\ep}_1 [j,k]=m^{\ep}_1 (t_1, t_2)$ and $m^{\ep}_3 [j,k] = m^{\ep}_3 (t_1, t_2)$ with $t_1 = 2^{-j} (k-1)$ and $t_2 = 2^{-j} k$ in relation \eqref{ineq:epbeta-mep-bounds}. The instants $t_1, t_2$ considered in \eqref{ineq:xi-ep-i-norm} are such that $t_2-t_1 = 2^{-j}$. 
Hence, there exists a constant $C$ such that
\begin{equation}\label{ineq:xi-ep-i-norm-ver}
\norm{\xi_{\ep , i}}_Z \leq C K \sum_{j=i}^{v}  j \sqrt{2^{-j}} \leq C K 2^{-i/2} 
\sum_{j=i}^{\infty} j 2^{(i-j)/2} \leq C i 2^{-i/2}.
\end{equation}
The quantity on the right in the above can be made arbitrarily small as 
$i \to \infty$ so that given any rational $\theta, \eta >0$, we can have 
$C i2^{-i/2} \leq \theta / \ln(1+\eta^{-1})$.
Therefore, there exists $\ep_0 (\omega) >0$ such that for all $\ep \leq \ep_0 (\omega)$, we can write~\eqref{ineq:xi-ep-i-norm-ver} as
\begin{equation*}
\norm{\xi_{\ep , i}}_Z \leq C i2^{-i/2} \leq \theta / \ln(1+\eta^{-1}) \,.
\end{equation*}

At the beginning of the proof, we have defined $v = v(\ep)$ such that $\lim_{\ep \to 0} v (\ep) = \infty$. Now let us recall the definition \eqref{eqdef:xi-ep-i} of $\xi_{\ep , i}$ and set $\delta = 2^{-i}$ for $i \leq v$ large enough. Invoking Markov’s inequality and choosing $c = \inf \{c : \E_Z [\phi (\xi_{\ep , i} /c)] \leq 1\}$ we get
\begin{eqnarray}
\PP_Z \left[\bigvee_{s, t \in F_v , |t-s|\leq \delta} |G^{\ep} (t) - G^{\ep} (s)| > \theta \right] 
&=& \PP_Z \left[\xi_{\ep , i} > \theta\right] = \PP_Z \left[\phi (\xi_{\ep , i} /c) > \phi (\theta /c)\right] \notag \\
&\leq& \frac{\E_Z [\phi (\xi_{\ep , i} /c)]}{\phi (\theta /c)} \leq
\frac{1}{\phi (\theta / \norm{\xi_{\ep , i}}_Z)} = \frac{1}{e^{\theta / \norm{\xi_{\ep , i}}_Z} -1} \notag \\
&\leq& \frac{1}{\exp \{\theta \,\cdot\, \ln (1+\eta^{-1}) / \theta\} -1} \notag \\
&=& \eta \,. \notag
\end{eqnarray}
This proves the desired inequality \eqref{ineq:prob-sup-diff-G-ep}, which is our claim.
\end{proof}

Now we proceed to prove \eqref{lim:max-increment-G-ep-v} 
which will complete the proof of our quenched tightness. Note that our method of proof here is different in nature from \cite{RR}. It relies mostly on martingale techniques together with a convenient domain decomposition (see \cref{fig:A_i-region-add}) for the point process $M^\ep$.

\begin{lemma}[Proof of \eqref{lim:max-increment-G-ep-v}]
\label{lemma:tightness-ineqstate2}
Let $\{G^{\ep} \,;\, \ep >0\}$ be the family of processes such that $G^{\ep} (t)$ is defined by \eqref{centered-Nep}.
Then as $\ep \to 0$ and $v = v(\ep) \to \infty$, the limit \eqref{lim:max-increment-G-ep-v} holds true in $\PP_Z$-probability, namely,
\begin{equation*}
\bigvee_{k =0}^{2^v-1} \sup_{0 \leq u \leq 2^{-v}} |G^{\ep} (k2^{-v} +u) - G^{\ep} (k 2^{-v})| \to 0 \,.
\end{equation*}
\end{lemma}

\begin{figure}
\begin{tikzpicture}[xscale = 1.6]
\path [fill = lime!60] (1.7,0) -- (1.7,2.3) -- (0,4) -- (0,1.7);
\path [fill = cyan!70] (1.7,2.3) -- (1.7,4.8) -- (0,4.8) -- (0,4);
\path [fill = orange!60] (4,0) -- (4,4.8) -- (1.7,4.8) -- (1.7,2.3);
\path [fill = magenta!40] (4,0) -- (1.7,2.3) -- (1.7,0);

\draw [very thick] [<->] (0,4.9) -- (0,0) -- (5.8,0);
\draw [thick] (0,1.7) -- (1.7,0);
\node [left] at (0,1.7) {$t_1$};
\draw [thick] (0, 4) -- (4 ,0);
\node [left] at (0,4) {$t_2$};
\draw [thick] (1.7,0) -- (1.7, 4.8);
\node [below] at (1.7,0) {$t_1$};
\draw [thick] (4, 0) -- (4,4.8);
\node [below] at (4,0) {$t_2$};
\node [above left] at (0,4.9) {$l$};
\node [below right] at (5.8,0) {$\gamma$};


\node at (0.75, 1.92) {$A_1$};
\node at (1, 4) {$A_2$};
\node at (3, 3) {$A_3$};
\node at (2.4, 0.8) {$A_5$};
\node at (5.4, 2) {$A_4 = A_3\cup A_5$};
\draw (3.8, 1.2) to (4.5, 1.8);
\end{tikzpicture}
\caption{Decomposition of $\R_+ \times \R_+$.}
\label{fig:A_i-region-add}
\end{figure}

\begin{proof}

Using the notation $m^{\ep}_i = m^{\ep}_i (t_1, t_2)$ and 
$A_i = A_i (t_1,t_2)$ introduced for relation~\eqref{eq:ep-m-1,3-int}
with $t_1 = k2^{-v}$ and $t_2 = k2^{-v}+u$,
similar to \eqref{conc.G.diff}, we have
\begin{eqnarray}
&&G^{\ep} (k2^{-v} +u) - G^{\ep} (k 2^{-v}) \nonumber\\
&=&
\widetilde{M^\ep}(A_3(k2^{-v}, k2^{-v}+u))-
\widetilde{M^\ep}(A_1(k2^{-v}, k2^{-v}+u))\nonumber\\
&=&
\ep^{\beta/2} 
\Big(M^{\ep}[ A_3 (k2^{-v}, k2^{-v} +u) ]-\E_Z \big[M^{\ep}[ A_3 (k2^{-v}, k2^{-v} +u) ] \big]\Big)\nonumber\\
&&-\ep^{\beta/2} 
\Big(M^{\ep}[ A_1 (k2^{-v}, k2^{-v} +u) ]-\E_Z \big[M^{\ep}[ A_1 (k2^{-v}, k2^{-v} +u) ] \big]\Big).
\label{eq:diff-G-ep-expand1s}
\end{eqnarray}
Now we analyze the individual summands of the above.
Note that the set $A_1(t_1, t_2)$ is {\em increasing} as a function of $t_2$, 
or $u$; see Figure \ref{fig:A_i-region-add}. More precisely,
$A_1 (t_1, t_2)\subset A_1 (t_1, t_2')$ for $t_2 < t_2'$.
Hence, due to
the independence property of Poisson point process, the following
process
$$
u \mapsto
{\mathfrak M}_1^\ep(u):=M^{\ep}[ A_1 (k2^{-v}, k2^{-v} +u) ]-\E_Z \big[M^{\ep}[ A_1 (k2^{-v}, k2^{-v} +u) ] \big]
$$
is a martingale, i.e.
$ \E_Z[{\mathfrak M}_1^\ep(u_2)\big|\mathcal F_{u_1}]=
{\mathfrak M}_1^\ep(u_1)
$ for all $u_1 < u_2$ where $\mathcal F_{u}$ is the $\sigma$-field generated
by the region $A_1(k2^{-v},k2^{-v}+u)$.

For the set $A_3(t_1,t_2)$, we can write it as
$A_3(t_1,t_2)=A_4(t_1,t_2)\backslash A_5(t_1,t_2)$,
where
\begin{eqnarray*}
A_4(t_1,t_2) &:=& \left\{ (\ga, l) \in \R_+ \times \R_+ : \, t_1 \leq \ga \leq t_2 \,\,\text{and}\,\, l\geq 0\right\},\\
A_5(t_1,t_2) &:=& \left\{(\ga, l) \in \R_+ \times \R_+ : \, t_1 \leq \ga \leq t_2 \,\,\text{and}\,\, \ga + l \leq t_2\right\}.
\end{eqnarray*}
Note that the sets $A_4(t_1,t_2)$ and $A_5(t_1,t_2)$ are again increasing in 
$t_2$, or $u$. Hence for the same reason,
the following is also a martingale for $i=4, 5$:
$$
u \mapsto
{\mathfrak M}_i^\ep(u):=M^{\ep}[ A_i (k2^{-v}, k2^{-v} +u) ]-\E_Z \big[M^{\ep}[ A_i (k2^{-v}, k2^{-v} +u) ] \big]
$$
(But note the subtle point that the filtrations for $\mathfrak M_4$ and
$\mathfrak M_5$ are generated by the regions 
$A_4(k2^{-v},k2^{-v}+u)$
and $A_5(k2^{-v},k2^{-v}+u)$ respectively, and hence are different).

With the above definition of sets, we can rewrite \eqref{eq:diff-G-ep-expand1s} as
\begin{equation}\label{defeq:mathfrak-M-increment}
G^{\ep} (k2^{-v} +u) - G^{\ep} (k 2^{-v})
= \ep^{\beta/2}\Big(
{\mathfrak M}_4^\ep(u) -{\mathfrak M}_5^\ep(u) -{\mathfrak M}_1^\ep(u)
\Big).
\end{equation}
Invoking \cite[Proposition 3.2]{DImparato}, we state that
for $i=1,4,5$ and with a constant $C$,
\begin{eqnarray}
\norm{\sup_{0 \leq u \leq 2^{-v}}
\left|\mathfrak{M}^\ep_i(u)\right| }_Z
\leq C
\norm{\mathfrak{M}^{\ep}_i(2^{-v})}_Z.
\label{ineq:sup-M-ep-i(u)}
\end{eqnarray}
Then we apply \cref{lemma:Poisson-rv-Orlicz,lemma:unifbound-epbeta-m-i} 
to $\mathfrak{M}^{\ep}_i (2^{-v})$. This yields,
for all $\ep \leq \ep_0 (\omega)$ small enough, 
$K$ a constant (as in \cref{lemma:Poisson-rv-Orlicz}), and $C$ a constant (as in \cref{lemma:unifbound-epbeta-m-i}),
\begin{eqnarray}
\norm{\mathfrak{M}^\ep_i(2^{-v})}_Z & = & 
\norm{
\left|M^{\ep}[ A_i (k2^{-v}, (k+1)2^{-v}) ]
-\E_Z \big[M^{\ep}[ A_i (k2^{-v}, (k+1)2^{-v}) ] \big]
\right|
}_Z\nonumber\\
&\leq &
K\sqrt{m^\ep_i(t_1,t_2)}
\leq
K \sqrt{C \ep^{-\beta} 2^{-v}} =
K \, \ep^{-\beta /2} 2^{-v/2} \,.
\label{ineq:bdd-M-ep-i-new}
\end{eqnarray}
(The statement of \cref{lemma:unifbound-epbeta-m-i}, stated for $i=1,3$ only, 
can be easily extended to $i=4, 5$.)
Plugging \eqref{ineq:bdd-M-ep-i-new} in \eqref{ineq:sup-M-ep-i(u)} 
and then in \eqref{defeq:mathfrak-M-increment}, we get
\begin{equation*}
\norm{\sup_{0 \leq u \leq 2^{-v}} \left| G^{\ep} (k2^{-v} +u) - G^{\ep} (k 2^{-v}) \right| }_Z \leq \ep^{\beta/2} K \, \ep^{-\beta /2} 2^{-v/2} = K \, 2^{-v/2} \,.
\end{equation*}

\noindent
In order to achieve \eqref{lim:max-increment-G-ep-v}, we still have to handle the sup operation $\bigvee_{k =0}^{2^v-1}$. However, this is done exactly as in the proof of \cref{lemma:m-j-Orlicz-norm-bounds}. 
Specifically, with \cite[Lemma 2.2.2, page 96]{VW} under the quenched probability $\PP_Z$, we have that for a constant $C_{\phi}$ only depending on $\phi$,
\begin{equation*}
\norm{\bigvee_{k =0}^{2^v-1} \sup_{0 \leq u \leq 2^{-v}} \left| G^{\ep} (k2^{-v} +u) - G^{\ep} (k 2^{-v}) \right| }_Z
\leq C_{\phi} \, \phi^{-1}(2^v -1) K \, 2^{-v/2} \,.
\end{equation*}
Recall the fact that $\phi^{-1} (x) = \ln (x+1)$ for $x \geq 0$. We also invoke an elementary fact that $\lim_{v \to \infty} v 2^{-v/2} = 0$. Upon rewriting the above, we have as $\ep \to 0$,
\begin{equation*}
\norm{\bigvee_{k =0}^{2^v-1} \sup_{0 \leq u \leq 2^{-v}} \left| G^{\ep} (k2^{-v} +u) - G^{\ep} (k 2^{-v}) \right| }_Z
\leq C_{\phi} \, v \ln 2 \, K \, 2^{-v/2} \to 0\,.
\end{equation*}
This verifies our claim \eqref{lim:max-increment-G-ep-v}.
\end{proof}

\section{Proof of \cref{prop:tightness-annclt}:\\ Tightness in the Annealed Setting}~\label{app:annealed-tight}
The proof largely follows \cref{app:quenched-tight} and~\cite[Section 3.1]{RR}, albeit under the annealed measure $\PP := \mathbf P \otimes \PP_Z$. Recall the definition of the Orlicz norm under the annealed measure in \cref{defnthm:orlicz-norm}. Paralleling~\cite[Lemma 1]{RR} we have 

\begin{lemma}~\label{lemma:non-Poisson-rv-Orlicz}
	Let $N^\ep(t)$ be the number of jobs in the system defined by~\eqref{eq:N-t1} with quenched mean $m^\ep (t)$ for a given strictly positive $\beta$. We assume that \cref{hyp:arr-times-Gamma,hyp:F-bar-increments,hyp:natual-annealed-bound,supercri-hyp:coef-la-psi} are satisfied and consider $t \geq t_0 >0$. Then, there exists $\ep_0 > 0$ and a universal deterministic constant $C$ such that for all $\ep \leq \ep_0$ we have	
	\begin{equation}\label{ineq:norm-C-mep-lemmac}
\|N^\ep (t) - m^\ep (t)\| \leq C \ep^{-\beta/2} \,.
	\end{equation}
\end{lemma}

\begin{remark}
As mentioned above, \cref{lemma:non-Poisson-rv-Orlicz} is inspired by \cite[Lemma 1]{RR}. However, \cite[Lemma 1]{RR} only applies to the case where $N^\ep(t)$ is a Poisson random variable under the annealed measure. This is clearly not the case in our random parameter setting. Therefore, our claim necessitates a separate proof for \cref{lemma:non-Poisson-rv-Orlicz}, by using the corresponding quenched result in \cref{lemma:Poisson-rv-Orlicz}. 
In particular, the expressions \eqref{eq:const-orlicz-Poisson1} and \eqref{eq:const-orlicz-Poisson2} will be implicitly used in the annealed analysis where the $Z$ dependence in the constants $K_{Z,1}$ and $K_{Z,2}$ are
eliminated.
\end{remark}

\begin{proof}[Proof of \cref{lemma:non-Poisson-rv-Orlicz}]
Let $C$ be a generic parameter. As mentioned in \cref{defnthm:orlicz-norm}, we wish to find the smallest $C$ such that
$$
\E \left[ e^{\frac{|N^\ep (t) - m^\ep (t)|}{C}}\right] \leq 2 \,.
$$
We now use the quenched computations and write
\begin{equation}\label{eq:rewrite-norm-lemmaC1}
\E \left[ e^{\frac{|N^\ep (t) - m^\ep (t)|}{C}}\right] = \E \left\{\E_Z \left[ e^{\frac{|N^\ep (t) - m^\ep (t)|}{C}}\right] \right\}.
\end{equation}
To handle the right hand side of \eqref{eq:rewrite-norm-lemmaC1}, we need two ingredients:

\noindent
\emph{(i)} In order to apply \cref{lemma:Poisson-rv-Orlicz}, we must have $m^\ep (t) \geq 1$. Therefore, we proceed to show that
there exists $\ep_0$ such that $m^\ep (t) \geq 1$ for all $\ep \leq \ep_0$. Now recall from \eqref{def:m-ep-expansionwith-m-ep-0} that
\begin{equation}\label{def:superm-ep-expansionwith-m-ep-0}
m^{\ep} (t) = \frac{1}{\ep^\beta}\, m^{\ep}_0 (t) \,,
\quad \text{with} \quad m^{\ep}_0 (t) = \iot h (s,t)\, \psi (Z_{s/\ep}) \,\dd s \,.
\end{equation}
Moreover, according to \cref{supercri-hyp:coef-la-psi} we have $\psi (Z_{s/\ep}) \geq \underline\psi$, and relation \eqref{def:F-bar-generalization} yields
\begin{equation*}
h (s,t) = \la (s) \bar F_s (t-s) \geq \frac{C}{(1+ (t-s))^\al} \,, \quad \text{with some constant} \,\, C\,.
\end{equation*}
Reporting this information in \eqref{def:superm-ep-expansionwith-m-ep-0}, it is easily seen that if $t \geq t_0$ we have
\begin{equation}\label{ineq:m-ep-0lemmaC1}
m^\ep_0 (t) \geq C \iot \frac{\dd s}{(1+s)^\al} \geq C_0 >0 \,,\quad
\text{where} \quad C_0 = C \int^{t_0}_0 \frac{\dd s}{(1+s)^\al} \,.
\end{equation}
Going back to \eqref{def:superm-ep-expansionwith-m-ep-0}, it is readily checked that if $t \geq t_0$, then $m^\ep (t) \geq C_0 / \ep^\beta$. Therefore, we deduce the existence of $\ep_0 >0$ such that $m^\ep (t) \geq 1$ for all $\ep \leq \ep_0$.
Owing to \cref{lemma:Poisson-rv-Orlicz}, we thus get that for $\ep < \ep_0$,
\begin{equation}\label{ineq:E_Z-exp-frac-lemmac1}
\E_Z \left[ \exp \left( \frac{|N^\ep (t) - m^\ep (t)|}{K \sqrt{m^\ep (t)}}\right)\right] \leq 2 \,.
\end{equation}
Let us highlight the fact that $\ep_0$ above is a deterministic constant.

\noindent
\emph{(ii)} Starting from \eqref{def:superm-ep-expansionwith-m-ep-0}, it is also easily checked that the following upper bound counterpart of \eqref{ineq:m-ep-0lemmaC1} holds true:
$$m^\ep (t) \leq \frac{C_1}{\ep^\beta}\,.$$
Hence, it stems from \eqref{ineq:E_Z-exp-frac-lemmac1} that
\begin{equation}\label{ineq:E_Z-exp-frac2lemmac1}
\E_Z \left[\exp \left(\frac{|N^\ep (t) - m^\ep (t)|}{C_2 \,\ep^{-\beta/2}}\right)\right] \leq 2\,,
\end{equation}
where $C_2 = K \sqrt{C_1}$. Plugging \eqref{ineq:E_Z-exp-frac2lemmac1} into \eqref{eq:rewrite-norm-lemmaC1}, our conclusion \eqref{ineq:norm-C-mep-lemmac} follows.
\end{proof}

Now we proceed to prove the tightness in the annealed setting.
\begin{proof}[Proof of \cref{prop:tightness-annclt}]
We first observe that tightness on $D[0, \infty)$ stems from tightness on any finite interval of the form $D[a,b]$, with $0 <a<b<\infty$. Therefore, it is sufficient to show that for any $\ep > 0, ~\eta > 0$ and $0 <a<b<\infty$, there exist $\ep_0 ,\delta > 0$ such that 
\begin{equation*}
	\PP \left( \underset{|s-t| < \delta,~ a\leq s,\,t \leq b}{\sup} |G^\ep(t) - G^\ep(s)| > \ep \right) < \eta,~\text{for all}~\ep \leq \ep_0 \,.
\end{equation*}
The proof of this claim follows directly from the verification of equations (3.13) and (3.14) in~\cite{RR}. 
With  \cref{lemma:non-Poisson-rv-Orlicz} in hand, this is a mere elaboration of
\cite[Section 3.1]{RR}. Details are left to the reader for the sake of conciseness.
\end{proof}

\bigskip

\bibliographystyle{plain}
\bibliography{references,references-1.bib}

\begin{thebibliography}{10}

\bibitem{anderson2016functional}
David Anderson, Joke Blom, Michel Mandjes, Halldora Thorsdottir, and Koen
  De~Turck.
\newblock A functional central limit theorem for a {Markov}-modulated
  infinite-server queue.
\newblock {\em Methodology and Computing in Applied Probability},
  18(1):153--168, 2016.

\bibitem{ArMoRa2005}
G{\'e}rard~Ben Arous, Stanislav Molchanov, and Alejandro~F Ram{\'\i}rez.
\newblock Transition from the annealed to the quenched asymptotics for a random
  walk on random obstacles.
\newblock {\em The Annals of Probability}, 33(6):2149--2187, 2005.

\bibitem{Bill}
Patrick Billingsley.
\newblock {\em Convergence of probability measures}.
\newblock John Wiley \& Sons, 1968.

\bibitem{blanchet2013continuous}
Jose Blanchet and Xinyun Chen.
\newblock Continuous-time modeling of bid-ask spread and price dynamics in
  limit order books.
\newblock {\em arXiv preprint arXiv:1310.1103}, 2013.

\bibitem{blom2014markov}
J~Blom, Offer Kella, Michel Mandjes, and Halldora Thorsdottir.
\newblock {Markov}-modulated infinite-server queues with general service times.
\newblock {\em Queueing Systems}, 76(4):403--424, 2014.

\bibitem{blom2013time}
J~Blom, M~Mandjes, and H~Thorsdottir.
\newblock Time-scaling limits for {Markov}-modulated infinite-server queues.
\newblock {\em Stochastic Models}, 29(1):112--127, 2013.

\bibitem{blom2016functional}
Joke Blom, Koen De~Turck, and Michel Mandjes.
\newblock Functional central limit theorems for {Markov}-modulated
  infinite-server systems.
\newblock {\em Mathematical Methods of Operations Research}, 83(3):351--372,
  2016.

\bibitem{blom2013large}
Joke Blom and Michel Mandjes.
\newblock A large-deviations analysis of {Markov}-modulated infinite-server
  queues.
\newblock {\em Operations Research Letters}, 41(3):220--225, 2013.

\bibitem{boxma2019infinite}
Onno Boxma, Offer Kella, and Michel Mandjes.
\newblock Infinite-server systems with {Coxian} arrivals.
\newblock {\em Queueing Systems}, 92(3-4):233--255, 2019.

\bibitem{CaNua}
Simon Campese, Ivan Nourdin, and David Nualart.
\newblock Continuous {Breuer-Major} theorem: Tightness and nonstationarity.
\newblock {\em The Annals of Probability}, 48(1):147--177, 2020.

\bibitem{CCG}
Patrick Cattiaux, Djalil Chafa\"{i}, and Arnaud Guillin.
\newblock {Central limit theorems for additive functionals of ergodic Markov
  diffusions processes}.
\newblock {\em ALEA: Latin American Journal of Probability and Mathematical
  Statistics}, 9(2):337--382, 2012.

\bibitem{Cin}
Erhan \c{C}inlar.
\newblock {\em Probability and Stochastics}.
\newblock Springer, 2011.

\bibitem{CH1}
Prakash Chakraborty and Harsha Honnappa.
\newblock A many-server functional strong law for a non-stationary loss model.
\newblock {\em Submitted}, 2019.

\bibitem{CMRW}
GL~Choudhury, A~Mandelbaum, MI~Reiman, and W~Whitt.
\newblock Fluid and diffusion limits for queues in slowly changing
  environments.
\newblock {\em Stochastic Models}, 13(1):121--146, 1997.

\bibitem{coffman1995polling}
EG~Coffman~Jr, AA~Puhalskii, and MI~Reiman.
\newblock Polling systems with zero switchover times: a heavy-traffic averaging
  principle.
\newblock {\em The Annals of Applied Probability}, pages 681--719, 1995.

\bibitem{de2014large}
KEES De~Turck and MRH Mandjes.
\newblock Large deviations of an infinite-server system with a linearly scaled
  background process.
\newblock {\em Performance Evaluation}, 75:36--49, 2014.

\bibitem{dean2020functional}
Justin Dean, Ayalvadi Ganesh, and Edward Crane.
\newblock Functional large deviations for {Cox} processes and {Cox$/G/\infty$}
  queues, with a biological application.
\newblock {\em Annals of Applied Probability}, 30(5):2465--2490, 2020.

\bibitem{DGQuenched}
Dmitry Dolgopyat and Ilya Goldsheid.
\newblock {Quenched limit theorems for nearest neighbour random walks in 1D
  random environment}.
\newblock {\em Communications in Mathematical Physics}, 315(1):241--277, 2012.

\bibitem{subcritical}
Stephen~G Eick, William~A Massey, and Ward Whitt.
\newblock The physics of the {$M_t/G/\infty$} queue.
\newblock {\em Operations Research}, 41(4):731--742, 1993.

\bibitem{fibich2012averaging}
Gadi Fibich, Arieh Gavious, and Eilon Solan.
\newblock Averaging principle for second-order approximation of heterogeneous
  models with homogeneous models.
\newblock {\em Proceedings of the National Academy of Sciences},
  109(48):19545--19550, 2012.

\bibitem{fouque2018optimal}
Jean-Pierre Fouque and Ruimeng Hu.
\newblock Optimal portfolio under fast mean-reverting fractional stochastic
  environment.
\newblock {\em SIAM Journal on Financial Mathematics}, 9(2):564--601, 2018.

\bibitem{fouque1999financial}
Jean-Pierre Fouque, George Papanicolaou, and K~Ronnie Sircar.
\newblock Financial modeling in a fast mean-reverting stochastic volatility
  environment.
\newblock {\em Asia-Pacific Financial Markets}, 6(1):37--48, 1999.

\bibitem{fralix2009infinite}
Brian~H. Fralix and Ivo~J.B.F. Adan.
\newblock An infinite-server queue influenced by a semi-{Markovian}
  environment.
\newblock {\em Queueing Systems}, 61(1):65--84, 2009.

\bibitem{HM}
Mariska Heemskerk and Michel Mandjes.
\newblock Exact asymptotics in an infinite-server system with overdispersed
  input.
\newblock {\em Operations Research Letters}, 47(6):513--520, 2019.

\bibitem{HvLM}
Mariska Heemskerk, Johan van Leeuwaarden, and Michel Mandjes.
\newblock Scaling limits for infinite-server systems in a random environment.
\newblock {\em Stochastic Systems}, 7(1):1--31, 2017.

\bibitem{hellings2012semi}
Ton Hellings, Michel Mandjes, and Joke Blom.
\newblock Semi-{Markov}-modulated infinite-server queues: approximations by
  time-scaling.
\newblock {\em Stochastic Models}, 28(3):452--477, 2012.

\bibitem{honnappa}
Harsha Honnappa, Rahul Jain, and Amy~R Ward.
\newblock A queueing model with independent arrivals, and its fluid and
  diffusion limits.
\newblock {\em Queueing Systems}, 80(1-2):71--103, 2015.

\bibitem{hunt1994large}
PJ~Hunt and TG~Kurtz.
\newblock Large loss networks.
\newblock {\em Stochastic Processes and their Applications}, 53(2):363--378,
  1994.

\bibitem{DImparato}
Daniele Imparato.
\newblock Martingale inequalities in exponential {Orlicz} spaces.
\newblock {\em Journal of Inequalities in Pure and Applied Mathematics}, 10,
  2009.

\bibitem{jansen2016large}
Hermanus~Marinus Jansen, MRH Mandjes, Koen De~Turck, and Sabine Wittevrongel.
\newblock A large deviations principle for infinite-server queues in a random
  environment.
\newblock {\em Queueing Systems}, 82(1-2):199--235, 2016.

\bibitem{jansen2019diffusion}
HM~Jansen, Michel Mandjes, Koen De~Turck, and Sabine Wittevrongel.
\newblock Diffusion limits for networks of markov-modulated infinite-server
  queues.
\newblock {\em Performance Evaluation}, 135:102039, 2019.

\bibitem{KYZ1}
RZ~Khashinskii, G~Yin, and Qing Zhang.
\newblock Asymptotic expansions of singularly perturbed systems involving
  rapidly fluctuating markov chains.
\newblock {\em SIAM Journal on Applied Mathematics}, 56(1):277--293, 1996.

\bibitem{khasminskii2004averaging}
Rafail~Z Khasminskii and George Yin.
\newblock On averaging principles: An asymptotic expansion approach.
\newblock {\em SIAM Journal on Mathematical Analysis}, 35(6):1534--1560, 2004.

\bibitem{khasminskij1968principle}
RZ~Khasminskii.
\newblock On the principle of averaging the {Ito's} stochastic differential
  equations.
\newblock {\em Kybernetika}, 4(3):260--279, 1968.

\bibitem{KYZ2}
RZ~Khasminskii, G~Yin, and Q~Zhang.
\newblock Constructing asymptotic series for probability distributions of
  markov chains with weak and strong interactions.
\newblock {\em Quarterly of Applied Mathematics}, 55(1):177--200, 1997.

\bibitem{KW1}
Song-Hee Kim and Ward Whitt.
\newblock Are call center and hospital arrivals well modeled by nonhomogeneous
  poisson processes?
\newblock {\em Manufacturing \& Service Operations Management}, 16(3):464--480,
  2014.

\bibitem{KV}
C.~Kipnis and S.~Varadhan.
\newblock Central limit theorem for additive functionals of reversible markov
  processes and applications to simple exclusions.
\newblock {\em Communications in Mathematical Physics}, 104:1--19, 1986.

\bibitem{KLO}
Tomasz Komorowski, Claudio Landim, and Stefano Olla.
\newblock {\em Fluctuations in Markov Processes: Time Symmetry and Martingale
  Approximation}.
\newblock Springer, 2012.

\bibitem{koops2017networks}
David~T Koops, Onno~J Boxma, and MRH Mandjes.
\newblock {Networks of $\cdot/G/\infty$ queues with shot-noise-driven arrival
  intensities}.
\newblock {\em Queueing Systems}, 86(3):301--325, 2017.

\bibitem{kurtz1992averaging}
Thomas~G Kurtz.
\newblock Averaging for martingale problems and stochastic approximation.
\newblock In {\em Applied Stochastic Analysis}, pages 186--209. Springer, 1992.

\bibitem{HLTY}
Yiran Liu, Harsha Honnappa, Samy Tindel, and Nung~Kwan Yip.
\newblock Infinite server queues in a random fast oscillatory environment.
\newblock {\em Queueing Systems}, 98:145--179, 2021.

\bibitem{lu2}
Yunan Liu and Ward Whitt.
\newblock {A many-server fluid limit for the $G_t/GI/s_t+ GI$ queueing model
  experiencing periods of overloading}.
\newblock {\em Operations Research Letters}, 40(5):307--312, 2012.

\bibitem{lu1}
Yunan Liu and Ward Whitt.
\newblock Many-server heavy-traffic limit for queues with time-varying
  parameters.
\newblock {\em The Annals of Applied Probability}, 24(1):378--421, 2014.

\bibitem{MM}
Avi Mandelbaum and William~A Massey.
\newblock Strong approximations for time-dependent queues.
\newblock {\em Mathematics of Operations Research}, 20(1):33--64, 1995.

\bibitem{MMR}
Avi Mandelbaum, William~A Massey, and Martin~I Reiman.
\newblock Strong approximations for {Markovian} service networks.
\newblock {\em Queueing Systems}, 30(1-2):149--201, 1998.

\bibitem{mandjes2016markov}
Michel Mandjes and Koen De~Turck.
\newblock Markov-modulated infinite-server queues driven by a common background
  process.
\newblock {\em Stochastic Models}, 32(2):206--232, 2016.

\bibitem{MW2}
William~A Massey and Ward Whitt.
\newblock Uniform acceleration expansions for {Markov} chains with time-varying
  rates.
\newblock {\em Annals of Applied Probability}, pages 1130--1155, 1998.

\bibitem{OP}
CA~O'cinneide and P~Purdue.
\newblock {The $M/M/\infty$ queue in a random environment}.
\newblock {\em Journal of Applied Probability}, 23(1):175--184, 1986.

\bibitem{PW}
Guodong Pang and Ward Whitt.
\newblock Two-parameter heavy-traffic limits for infinite-server queues.
\newblock {\em Queueing Systems}, 65:325--364, 2010.

\bibitem{pender1}
Jamol Pender.
\newblock Nonstationary loss queues via cumulant moment approximations.
\newblock {\em Probability in the Engineering and Informational Sciences},
  29(1):27--49, 2015.

\bibitem{pender2}
Jamol Pender and Young~Myoung Ko.
\newblock Approximations for the queue length distributions of time-varying
  many-server queues.
\newblock {\em INFORMS Journal on Computing}, 29(4):688--704, 2017.

\bibitem{perry2011ode}
Ohad Perry and Ward Whitt.
\newblock An {ODE} for an overloaded {X} model involving a stochastic averaging
  principle.
\newblock {\em Stochastic Systems}, 1(1):59--108, 2011.

\bibitem{Peterson08}
Jonathon Peterson.
\newblock Limiting distributions and large deviations for random walks in
  random environments.
\newblock 2008.
\newblock PhD Thesis.

\bibitem{prekopa}
Andr{\'a}s Pr{\'e}kopa.
\newblock {On secondary processes generated by a random point distribution of
  Poisson type}.
\newblock {\em Ann. Univ. Sci. Budapest Sectio Math}, 1:153--170, 1958.

\bibitem{RR}
Sidney Resnick and Holger Rootz{\'e}n.
\newblock Self-similar communication models and very heavy tails.
\newblock {\em Annals of Applied Probability}, 10(3):753--778, 2000.

\bibitem{SS}
Peter Spreij and Jaap Storm.
\newblock Diffusion limits for a {Markov} modulated counting process.
\newblock {\em arXiv preprint arXiv:1801.03682}, 2018.

\bibitem{VW}
Aad~W. van~der Vaart and Jon~A. Wellner.
\newblock {\em Weak Convergence and Empirical Processes with Applications to
  Statistics}.
\newblock Springer, 1996.

\bibitem{whitt1}
Ward Whitt.
\newblock Time-varying queues.
\newblock {\em Queueing models and service management}, 1(2), 2018.

\bibitem{ZHG2}
Zeyu Zheng, Harsha Honnappa, and Peter~W Glynn.
\newblock Approximating performance measures for slowly changing non-stationary
  {Markov} chains.
\newblock {\em arXiv preprint arXiv:1805.01662}, 2018.

\bibitem{ZHG1}
Zeyu Zheng, Harsha Honnappa, and Peter~W Glynn.
\newblock Approximating systems fed by {Poisson} processes with rapidly
  changing arrival rates.
\newblock {\em Operations Research}, 2021.

\end{thebibliography}
\end{document}